\documentclass[12pt] {article}
\usepackage{amsmath,amsthm}
\usepackage{amssymb,latexsym}
\title{Quaternionic Monge-Amp\`ere equations.}
\date{}
\author{ Semyon Alesker
\\  { \normalsize Department of Mathematics, Tel Aviv University, Ramat Aviv}
 \\  { \normalsize 69978 Tel Aviv,
Israel }
\\ {\normalsize e-mail: semyon@post.tau.ac.il}}
\newcommand{\RR}{\mbox{\rm $~\vrule height6.5pt width0.5pt
depth0.3pt\!\!$R}}

\newcommand{\CC}{\mbox{\rm $~\vrule height6.5pt width0.5pt
depth0.3pt\!\!$C}}

\newcommand{\HH}{\mbox{\rm $~\vrule height6.5pt width0.5pt
depth0.3pt\!\!$H}}

\def\eps{\varepsilon}

\def\lam{\lambda}

\def\str{\longrightarrow}

\def\qed { Q.E.D. }

\def\hn { \hat \nu }

\def\rnn{ \RR _{\geq 0}}

\def\psh{plurisubharmonic }

\def\db{\frac{\partial}{\partial \bar q}}
\def\dq{\frac{\partial}{\partial  q}}

\def\dfq{\frac{\partial ^2 f}{\partial\bar q_i \partial q_j}}

\swapnumbers
\newtheorem{theorem}{Theorem}[subsection]

\newtheorem{lemma}[theorem]{Lemma}
\newtheorem{proposition}[theorem]{Proposition}
\newtheorem{claim}[theorem]{Claim}
\theoremstyle{definition}

\newtheorem{definition}[theorem]{Definition}
\newtheorem{remark}[theorem]{Remark}
\theoremstyle{proposition-definition}
\newtheorem{proposition-definition}[theorem]{Proposition-Definition}
\def\bfm{{\cal B} (\phi, f)}
\def\bo{\partial \Omega}
\def\clo{\bar \Omega}
\def\vz{v_{\zeta}}
\def\fss{F^{**}}
\def\ggo{\Gamma_0}
\def\pt{\partial ^2}
\def\oe{\omega(\eps)}

\def\bb{\partial B}
\def\bo{\partial \Omega}
\def\du{\frac{\partial ^2 u}{\partial\bar q_i \partial q_j}}
\def\due{\frac{\partial ^2 u_{\eps}}{\partial\bar q_i \partial q_j}}
\def\dun{\frac{\partial ^2 u_N}{\partial\bar q_i \partial q_j}}

\begin{document}
\maketitle
\begin{abstract}
The main result of this paper is the existence and uniqueness of
solution of the Dirichlet problem for quaternionic Monge-Amp\`ere
equations in quaternionic strictly  pseudoconvex bounded domains
in $\HH ^n$. We continue the study of the theory of
plurisubharmonic functions of quaternionic variables started by
the author at \cite{alesker}.
\end{abstract}

 \setcounter{section}{-1}
\section{Introduction.}
\setcounter{subsection}{1} This paper is a continuation of
author's previous paper \cite{alesker}. In \cite{alesker} we have
developed the necessary algebraic technique and we have introduced
and studied the class of plurisubharmonic functions of
quaternionic variables (this class was independently introduced
also by G. Henkin \cite{henkin}). The main result of the present
paper is the existence of a generalized solution of the Dirichlet
problem for quaternionic Monge-Amp\`ere equations in quaternionic
strictly pseudoconvex bounded domains in $\HH ^n$. The uniqueness
of solution was established in \cite{alesker}.

The versions of this result for real and complex Monge-Amp\`ere
equations were established in classical papers by A.D. Aleksandrov
\cite{aleksandrov} (the real case) and E. Bedford and B. Taylor
\cite{bedford-taylor} (the complex case). We prove also a result
on the regularity of solution in the Euclidean ball. For real
Monge-Amp\`ere equations this result was proved by L. Caffarelli,
L. Nirenberg, and J. Spruck \cite{CNS} for arbitrary strictly
convex bounded domains, and for complex Monge-Amp\`ere equations
by L. Caffarelli, J. Kohn, L. Nirenberg, and J. Spruck \cite{CKNS}
and N. Krylov \cite{krylov2} for arbitrary strictly pseudoconvex
bounded domains.

The real Monge-Amp\`ere equations appear in various geometric
problems such as the Minkowski problem (see A. Pogorelov
\cite{pogorelov}) The Dirichlet problem has received considerable
study. The {\itshape interior} regularity of the solution of the
Dirichlet problem was proved by A. Pogorelov, and the proof was
briefly described  in \cite{pogorelov1}-\cite{pogorelov3}. The
complete proof was published in \cite{pogorelov} and
\cite{cheng-yau1}, \cite{cheng-yau2}. In \cite{cheng-yau2} Cheng
and Yau gave a different proof of interior regularity; they have
also studied some related geometric problems on flat manifolds.
Another motivation of studying of the real Monge-Amp\'ere
equations is the Monge-Kantorovich problem on the measure
transportation (see e.g. \cite{evans} and
\cite{rachev-ruschendorf}).

The complex Monge-Ampr\`ere equations were studied in particular
in the connection to K\"ahler geometry. Good references are the
books by T. Aubin \cite{aubin}, A. Besse \cite{besse}, D. Joyce
\cite{joyce-book}.

There is a general philosophy (promoted especially by V. Arnold,
see e.g. \cite{arnold}) that some mathematical theories should
have three versions which are analogous to each other in certain
sense: real, complex, and quaternionic. However they should
reflect different phenomena. In this paper we present the theory
of plurisubharmonic functions of quaternionic variables and the
theory of quaternionic Monge-Amp\`ere equations whose study we
have started in \cite{alesker}. Their real and complex analogs are
well known (the real analog of the theory of plurisubharmonic
functions is the theory of convex functions).

We are going to formulate our main result more precisely and
recall the main notions from \cite{alesker}. Let $\HH$ denote the
(non-commutative) field of quaternions. Let $\HH ^n$ denote the
space of $n$-tuples of quaternions
 $(q_1, q_2, \dots , q_n)$. We consider $\HH ^n$ as right $\HH$-module
(we call it right $\HH$-vector space). An $n\times n$ quaternionic
matrix $A=(a_{ij})$ is called {\itshape hyperhermitian} if $A^*
=A$, i.e. $a_{ij}= \bar a_{ji}$ for all $i,\, j$.

In order to write the classical (real or complex) Monge-Amp\`ere
equations one has to use the notion of determinant of matrices. By
now there is no construction of determinant of matrices with
non-commuting (even quaternionic) entries which would have
{\itshape all} the properties of the usual (commutative)
determinant. The most general theory of non-commutative
determinants is due to Gelfand and Retakh (see
\cite{gelfand-retakh1}-\cite{gelfand-retakh3}, also
\cite{gelfand-retakh-wilson}, \cite{ggrw}). However it turns out
that on the class of quaternionic hyperhermitian matrices there is
a notion of the Moore determinant which has all the properties of
the usual determinant of complex (resp. real) hermitian (resp.
symmetric) matrices. Some of these properties are reviewed in
Section 1, and we refer to \cite{alesker} for further details and
references. Here we mention only that the Moore determinant
depends polynomially on the entries of a hyperhermitian matrix,
and the Moore determinant of any complex hermitian matrix $A$
considered as quaternionic hyperhermitian coincides with the usual
determinant of $A$. We denote the Moore determinant of $A$ by
$\det A$.

The quaternionic Monge-Amp\`ere equation is written in terms of
this determinant. We have also to recall the notion of
plurisubharmonic function of quaternionic variables and the
definition of quaternionic strictly pseudoconvex domain following
\cite{alesker}.
 Let $\Omega $ be a domain in $\HH ^n$.
\begin{definition}
A real valued function $u: \Omega \str \RR$ is called quaternionic
plurisubharmonic (psh) if it is upper semi-continuous and its
restriction to any right {\itshape quaternionic} line is
subharmonic.
\end{definition}
 Recall that upper semi-continuity means that
 $f(x_0)\geq \underset{x\str x_0}{\limsup f(x)}$ for any $x_0\in \Omega$.

 Also we will call a $C^2$-smooth function $u:\Omega\str \RR$ to
 be {\itshape strictly plurisubharmonic} if its restriction to any
 right quaternionic line is strictly harmonic (i.e. the Laplacian
 is strictly positive).
\begin{definition}
  An open bounded domain $\Omega\subset
\HH^n$ with a smooth boundary $\partial \Omega$ is called strictly
pseudoconvex if for every point $z_0\in \partial \Omega$ there
exists a neighborhood ${\cal O}$ and a smooth strictly psh
function $h$ on ${\cal O}$ such that $\Omega \cap {\cal O}= \{ h<0
\}$, $ h(z_0)=0$, and $\nabla h(z_0) \ne 0$.
\end{definition}

We will write a quaternion $q$ in the usual form
$$q= t+ x\cdot i +y\cdot j+ z\cdot k ,$$
where $t,\, x,\, y,\, z$ are real numbers, and $i,\, j,\, k$
satisfy the usual relations
$$i^2=j^2=k^2=-1, \, ij=-ji=k,\, jk=-kj=i, \, ki=-ik=j.$$

The Dirac-Weyl operator $\frac {\partial}{\partial \bar q}$ is
defined as follows. For any $\HH$-valued function $f$
$$\db f:=\frac{\partial f}{\partial  t}  +
i \frac{\partial f}{\partial  x} + j \frac{\partial f}{\partial
y} + k \frac{\partial f}{\partial  z}.$$

Let us also define the operator $\dq$:
$$\dq f:=\overline{ \db \bar f}=
\frac{\partial f}{\partial  t}  -
 \frac{\partial f}{\partial x}  i-
 \frac{\partial f}{\partial  y} j-
\frac{\partial f}{\partial  z}  k.$$

{\bf Remarks.} (a) The operator $\db$ is called sometimes the
Cauchy-Riemann-Moisil-Fueter operator since it was introduced by
Moisil in \cite{moisil} and used by Fueter \cite{fueter1},
\cite{fueter2} to define the notion of quaternionic analyticity.
For further results on quaternionic analyticity we refer e.g. to
\cite{brackx-delanghe-sommen}, \cite{palamodov},
 \cite{pertici}, \cite{sudbery}, and for applications to mathematical
physics to \cite{gursey-tze}.
   Another used name for this
 operator is Dirac-Weyl operator. But in fact it was used earlier
 by J.C. Maxwell  in \cite{maxwell}, vol. II, pp.570-576, where he
 has applied the quaternions to electromagnetism.

 (b) Note that $$\db= \frac{\partial}{\partial t}+ \nabla,$$
 where $\nabla= i\frac{\partial}{\partial x}+j\frac{\partial}{\partial
 y}+ k\frac{\partial}{\partial z}$. The operator $\nabla$ was
 first introduced by W.R. Hamilton in \cite{hamilton}.

 (c) In quaternionic analysis one considers a right version of the
 operators $\db$ and $ \dq$ which are denoted respectively by
 $\overset{\leftarrow}{\db}$ and $ \overset{\leftarrow}{\dq}$. The
 operator $ \overset{\leftarrow}{\dq}$ is related to
 $\overset{\leftarrow}{\db}$ by the same formula as $\dq$ is
 related to $\db$, and $\overset{\leftarrow}{\db}$ is defined as
 $$\overset{\leftarrow}{\db}f :=\frac{\partial f}{\partial  t}  +
 \frac{\partial f}{\partial  x}i +  \frac{\partial f}{\partial y}j
+  \frac{\partial f}{\partial  z}k .$$ For a real valued function
$f$ the derivatives $\frac{\partial ^2 f}{\partial q_i \partial
\bar q_j}$ and $\frac{\partial ^2 f}{\overset{\leftarrow}{\partial
q_i} \overset{\leftarrow}{\partial \bar q_j}}$ are quaternionic
conjugate to each other.

Now we can write the quaternionic Monge-Amp\`ere equation with
respect to $C^2$- smooth psh function $u$ on $\Omega$:
$$\det (\du) =f,$$
where $f$ is a given function. Note that the matrix $(\du)$ is
quaternionic hyperhermitian (since $u$ is real valued), $\det$
means the Moore determinant of this matrix. Note also that since
the function $u$ is psh,  the matrix $(\du)$ is non-negative
definite, and hence its Moore determinant is non-negative (the
notion of positive definiteness of a hyperhermitian matrix is
recalled in Section 1, Definition 1.1.4).

One of the main results of Section 2 of \cite{alesker} was the
definition of non-negative measure also denoted by $\det (\du)$
for any {\itshape continuous} psh function $u$ (which is not
necessarily smooth). That construction generalizes to the
quaternionic situation the well known constructions in the real
and complex cases due respectively to A.D. Aleksandrov
\cite{aleksandrov} and Chern-Levine-Nirenberg
\cite{chern-levine-nirenberg}. Now we can formulate the main
results of this paper.
\begin{theorem}
Let $\Omega \subset \HH ^n$ be a bounded quaternionic strictly pseudoconvex
 domain. Let $f\in C(\bar \Omega), \, f \geq 0$. Let $\phi \in C(\partial \Omega)$.
Then there exists unique function $u\in C(\bar \Omega)$ which is psh in $\Omega$
such that
$$\det (\du)=f \mbox{ in } \Omega,$$
$$u|_{\partial \Omega}\equiv \phi.$$
\end{theorem} Note that the uniqueness was proved in \cite{alesker}.
Theorem 0.1.3 claims existence of a solution in a generalized
sense (e.g. the function $u$ does not have to be smooth). It is of
interest to prove the regularity of solution $u$ under assumptions
of regularity of the initial data $f,\, \phi$. We can prove it
when the domain $\Omega$ is the Euclidean ball $B$ in $\HH^n$. In
Section 6 we prove the following result (called Theorem 6.0.1).
\begin{theorem}
Let $f\in C^{\infty}(\bar B), \, f>0$. Let $\phi \in
C^{\infty}(\partial B)$. There exists unique psh function $u\in
C^{\infty}(\bar B)$ which is a solution of the Dirichlet problem
$$\det (\frac{\partial ^2 u}{\partial\bar q_i \partial q_j})=f,$$
$$u|_{\partial B}= \phi.$$
\end{theorem}
The real version of this result was proved for arbitrary strictly
convex bounded domains in $\RR^n$ by Caffarelli, Nirenberg, and
Spruck \cite{CNS}. The complex version of it was proved for
arbitrary strictly pseudoconvex bounded domains in $\CC^n$ by
Caffarelli, Kohn, Nirenberg, and Spruck \cite{CKNS} and Krylov
\cite{krylov2}. Our method is a modification of the method of the
last paper \cite{CKNS}. Note also that in the case $n=1$ the
problem is reduced to the classical Dirichlet problem for the
Laplacian in $\RR^4$ (which is a linear problem); it was solved in
XIX century.

Let us make a few comments why the method of \cite{CKNS} can not
be generalized immediately to arbitrary strictly pseudoconvex
bounded domain in $\HH^n$. The main difficulty is that in the
complex case one uses the holomorphic transformations to make the
domain to be (locally) close to the Euclidean ball. In the
quaternionic situation it does not work. Indeed in the complex
case the class of diffeomorphisms of a domain which are either
holomorphic or anti-holomorphic can be characterized as the class
of diffeomorphisms preserving the class of psh functions. In the
quaternionic situation the class of diffeomorphisms preserving psh
functions is very small: all of them must be affine
transformations, more precisely modulo translations the
corresponding group is equal to $GL_n(\HH) Sp(1)$. The last fact
is proved in Subsection 2.2.

Let us make a few comments on the method of the proof of Theorem
0.1.3. It uses the solution of the Dirichlet problem in the unit
ball given by Theorem 0.1.4. The method to deduce  the general
case from this one follows the lines of the paper
\cite{bedford-taylor} by Bedford and Taylor. It was necessary to
generalize to the quaternionic situation  many results from the
usual (complex) theory of plurisubharmonic functions (this
investigation was started in \cite{alesker}). Sections 4 and 5 of
this paper follow very closely the complex case
\cite{bedford-taylor}.

The paper is organized as follows. In Section 1 we review the
necessary facts from the theory on non-commutative determinants;
the exposition follows \cite{alesker}. In Section 2 we review the
theory of plurisubharmonic functions of quaternionic variables as
it was  developed in \cite{alesker}. In Section 3 we construct the
matrix valued measure
 $\left( \frac{\pt u}{{\partial\bar q_i \partial q_j}}\right)$
for any (finite) psh function $u$ on $\Omega$. This construction
is a quaternionic version of the well known analogous construction
in the complex case (see \cite{lelong}, p.70). This construction
will be used in the proof of Theorem 0.1.3. In Section 4 we
construct for any finite psh function $u$ an operator $\Phi(u)$
which is essentially $(\det \pt u)^{\frac{1}{n}}$ (following
Section 5 of the paper \cite{bedford-taylor} by Bedford and
Taylor). This operator plays an important technical role in the
proof of Theorem 0.1.3. In Section 5 we establish several facts on
the envelopes of functions from the Perron-Bremermann families
following closely again the technique developed in the complex
case in \cite{bedford-taylor}. In Section 6 we prove Theorem 0.1.4
on the existence of $C^{\infty}$-regular solution of the Dirichlet
problem in the unit ball under appropriate assumptions on the
regularity of the initial data. In Section 7 we prove Theorem
0.1.3; the proof uses the results of all previous sections. In
Section 8 we discuss further the notion of strictly pseudoconvex
domain in the quaternionic space. Thus we introduce the
quaternionic analogue of the Levi form and consider some examples.
In Subsection 8.3 we state some open questions.

{\bf Acknowledgements.} We are grateful to P. Biran, M. Gromov, G.
Henkin, D. Kazhdan, L. Nirenberg, V. Palamodov, and L. Polterovich
for useful and stimulating discussions. We express our gratitude
to M. Sodin for numerous explanations of the classical theory of
subharmonic functions. We  thank J.C. Wood for useful
correspondences regarding harmonic morphisms.

\section{Background from non-commutative linear algebra.}
\setcounter{subsection}{1}
\setcounter{theorem}{0}
In this section we review some material on non-commutative determinants. More
 precisely we will recall some facts on the Dieudonn\'e and Moore determinants
of quaternionic matrices following \cite{alesker}.

The
Dieudonn\'e determinant of quaternionic matrices
 behaves exactly like the absolute value of
the usual determinant of real or complex matrices from all points
of view (algebraic and analytic). Let us denote by $M_n(\HH)$ the
set of all quaternionic $(n\times n)$- matrices. The Dieudonn\'e
determinant $D$ is defined on this set and takes values in
non-negative real numbers:
$$D:M_n(\HH) \str \rnn .$$ Then one has the following (known) results:

\begin{theorem}
(i) For any complex $(n \times n)$- matrix $X$ considered as
quaternionic matrix, the Dieudonn\'e determinant $D(X)$ is equal
to the absolute value of the usual determinant of $X$.

(ii) For any quaternionic matrix $X$
$$D(X)=D(X^t)=D(X^*),$$
where $X^t$ and $X^*$ denote
the transposed and quaternionic conjugate matrices of $X$  respectively.

(iii) $D(X\cdot Y) =D(X) D(Y).$
\end{theorem}

The following result is a weak version of the decomposition of the
determinant in row (column).
\begin{theorem}
Let $A= \left[ \begin{array}{ccc}
a_{11}&\dots&a_{1n}\\
\multicolumn{3}{c}{\dotfill} \\
a_{n1}&\dots&a_{nn}
\end{array}
\right] $ be a quaternionic matrix.
Then $$D(A) \leq \sum_{i=1}^n |a_{1i}| D(M_{1i}).$$
Similar inequalities hold for any other row or column.
\end{theorem}

(In this theorem $|a|$ denotes the absolute value of a quaternion $a$,
and $M_{pq}$ denotes the minor of the matrix $A$ obtained from it
by deleting the $p$-th row and $q$-th column).

In a sense, the Dieudonn\'e determinant provides the theory of
{\itshape absolute value } of determinant. However it is not
always sufficient and we loose many of the algebraic properties of
the usual determinant. The notion of Moore determinant provides
such a theory, but only on the class of quaternionic {\itshape
hyperhermitian } matrices. Remind that a square quaternionic
matrix $A$ is called hyperhermitian if its quaternionic conjugate
$A^*$ is equal to $A$. The Moore determinant
 denoted by $det$
is defined on the class of all hyperhermitian matrices and takes
real values. For the construction of the Moore determinant we
refer to \cite{alesker}, Subsection 1.1, where one can also find
the references to the original papers.
 The important advantage of the Moore determinant with respect to the Dieudonn\'e
determinant is that it depends polynomially on the entries of a
matrix; it has already all the
 algebraic and analytic properties of the usual determinant of real symmetric
and complex hermitian matrices. Let us state some of them.

\begin{theorem}
(i) The Moore determinant of any
complex hermitian matrix considered as quaternionic hyperhermitian matrix
is equal to its usual determinant.

(ii) For any hyperhermitian matrix $A$ and any quaternionic matrix $C$
$$det (C^*AC)= detA \cdot det(C^*C).$$
\end{theorem}

{\bf Examples.}

(a) Let $A =diag(\lam_1, \dots, \lam _n)$ be a diagonal matrix
with real $\lam _i$'s. Then $A$ is hyperhermitian and its Moore
determinant $detA= \prod _i \lam_i$.

(b)  A general hyperhermitian $2 \times 2$ matrix $A$ has the form
 $$ A=  \left[ \begin {array}{cc}
                     a&q\\
                \bar q&b\\
                \end{array} \right] ,$$
where $a,b \in \RR, \, q \in \HH$. Then its Moore determinant
is equal to
$det A =ab - q \bar q$.

Let us remind the definition of positive definiteness of
hyperhermitian quaternionic matrix following \cite{alesker}.
\begin{definition}
Let $A=(a_{ij})_{i,j=1}^n$ be a hyperhermitian quaternionic
matrix. $A$ is called \itshape{ non-negative definite} if for
every $n$-column of quaternions $\xi =(\xi_i)_{i=1}^n$ one has
$$\xi^* A \xi= \sum_{i,j}\bar {\xi_i} a_{ij} \xi_j\geq 0.$$
Similarly $A$ is called \itshape{positive definite} if the above
expression is strictly positive once $\xi \ne 0$.
\end{definition}

 In terms of the Moore determinant one can prove the
generalization of the classical Sylvester criterion of positive
definiteness of hyperhermitian matrices (Theorem 1.1.13 in
\cite{alesker}). In terms of the Moore determinant one can
introduce the notion of the mixed discriminant and to prove the
analogues of Alexandrov's inequalities for mixed discriminants
(Theorem 1.1.15 and Corollary 1.1.16 in \cite{alesker}).

The (well known) relation between the Dieudonn\'e and Moore determinants
is as follows: for any hyperhermitian matrix $X$
$$D(X)=|det X|.$$

Note that the Dieudonn\'e determinant was introduced originally by J. Dieudonn\'e
in  \cite{dieudonne} (see also \cite{artin} for his theory).
 It can be defined for arbitrary (non-commutative) field.
On more modern language this result can be formulated as a
computation of the $K_1$-group of a non-commutative field (see
e.g. \cite{rosenberg}). Note also that there is a more recent very
general theory of non-commutative determinants (or
quasideterminants) due to I. Gelfand and V. Retakh generalizing in
certain direction the theory of the Dieudonn\'e determinant and
many other known theories of non-commutative determinants. First
it was introduced in \cite{gelfand-retakh1}, see also
\cite{gelfand-retakh3}, \cite{ggrw} and references therein for
further developments and applications. In a recent paper
\cite{gelfand-retakh-wilson} the connection of the Moore and
Dieudonn\'e determinants of quaternionic matrices to the theory of
quasideterminants was made very explicit and well understood.

We would also like to mention a different direction of a
development of the {\itshape quaternionic} linear algebra started
by D. Joyce \cite{joyce} and applied by himself to hypercomplex
algebraic geometry. We refer also to D. Quillen's paper
\cite{quillen} for further investigations in that direction.
Another attempt to understand the quaternionic linear algebra from
the topological point of view was done by the author in
\cite{alesker2}.

\section{Review of the theory of psh functions of quaternionic variables.}
\subsection{Some results from \cite{alesker}.}
\setcounter{theorem}{0} We recall the basic facts from the theory
of plurisubharmonic (psh) functions of quaternionic variables
established by the author in \cite{alesker} (see Definition 0.1.1
in this paper). The operators $\db$ and $\dq$ were defined in the
introduction.

First one has the following simple fact (see Proposition 2.1.6 in \cite{alesker}).
\begin{proposition}
A real valued twice continuously
differentiable function $f$ on the domain $\Omega \subset \HH ^n$
is quaternionic
plurisubharmonic if and only if at every point $q \in \Omega$
the matrix $(\frac{\partial ^2 f}{\partial\bar q_i \partial q_j})(q)$
is non-negative definite.
\end{proposition}

Note that the matrix in the statement of proposition is
quaternionic hyperhermitian (since the function $f$ is real
valued). The more important thing is that in analogy to the real
and complex cases one can define for any continuous quaternionic
\psh function $f$ a non-negative measure $det(\frac{\partial ^2
f}{\partial\bar q_i
\partial q_j})$, where $det$ denotes the Moore determinant
(this measure is obviously defined for smooth $f$). One has the
following continuity result.

\begin{theorem}
Let $\{f_N\}$ be a sequence of continuous quaternionic \psh
functions in a domain $\Omega \subset \HH^n$. Assume that this
sequence converges uniformly on compact subsets to a function $f$.
Then $f$ is continuous quaternionic \psh function. Moreover the
sequence of measures $det(\frac{\partial ^2 f_N}{\partial\bar q_i
\partial q_j})$ weakly converges to the measure $det(\dfq)$.
\end{theorem}

The proofs of analogous results in real and complex cases can be
found in \cite{aubin}, where the exposition of this topic follows
the approach of Chern-Levine-Nirenberg \cite{chern-levine-nirenberg} and
  Rauch-Taylor \cite{rauch-taylor}. For the complex
case we refer also to the classical book by P. Lelong \cite{lelong}.

The next result is called the {\itshape minimum principle} (Theorem 2.2.1 in \cite{alesker}).
\begin{theorem}
 Let $\Omega$ be a bounded  open set
in $\HH^n$. Let $u,\, v$ be continuous functions on $\bar \Omega$
which are \psh in $\Omega$. Assume that
$$det(\frac{\partial ^2 u}{\partial\bar q_i \partial q_j})\leq
det(\frac{\partial ^2 v}{\partial\bar q_i \partial  q_j}) \mbox{
in } \Omega.$$ Then
$$min\{u(z)-v(z)|z\in \bar \Omega\} =min\{u(z)-v(z)|z\in \partial \Omega\}.$$
\end{theorem}
\subsection{Diffeomorphisms preserving psh functions.} In this
subsection we prove the following proposition.
\begin{proposition}
 Let $\Omega \subset \HH^n$ be  a domain. Let $F:\Omega \str
\Omega$ be a diffeomorphism such that for every open set ${\cal O}
\subset \Omega$ and for any psh function $f$ on ${\cal O}$ the
function $F^*f$ is psh on $F^{-1} ({\cal O})$. Then $F$ is an
affine transformation which can be written as a composition of a
translation and a linear transformation from the group
$GL_n(\HH)Sp(1)$.
\end{proposition}
In the statement of the theorem the group $GL_n(\HH)Sp(1)$ is
defined as follows. On the right quaternionic space $\HH^n$ there
is a left action of the  group of $\HH$-linear invertible
transformations $GL_n(\HH)$. Also the group $Sp(1)$ of norm one
quaternions acts on $\HH^n$ from the right. Both actions commute
and the group they generate is denoted by $GL_n(\HH)Sp(1)$. Note
that it is isomorphic to $(GL_n(\HH)\times Sp(1))/\{\pm Id\}$. Now
let us prove the proposition. Note also that all such affine
transformations preserving the domain $\Omega$ must preserve the
class of psh functions (see \cite{alesker}, Subsection 2.1).

{\bf Proof.} Let $U$ be any domain in $\HH^1$ and let $m:U\str
\Omega$ be any $\HH$-linear map. Let $p:\Omega \str \HH^1$ be any
$\HH$-linear projection. Consider the composition $p\circ F\circ
m:U\str \HH^1$. It is easy to see that this map preserves the
class of psh functions which is one-dimensional case means just
that it preserves the class of subharmonic functions. Hence
$p\circ F\circ m$ preserves the class of harmonic functions (i.e.
it is so called {\itshape harmonic morphism}, see e.g. \cite{wood}
for more details and references) . However there is a general
result of B. Fuglede \cite{fuglede} which says the following. Let
$g:M\str N$ be a smooth map between Riemannian manifolds of the
same dimension greater than 2 which preserves the class of
harmonic functions (in the above sense). Then $g$ is a conformal
mapping with the constant coefficient of conformality. When $M$
and $N$ are linear vector spaces with Euclidean metrics this
result together with the classical Liuville theorem imply that $g$
is a composition of homothety, translation, and orthogonal
transformation.

It easily follows that our original map $F$ is an affine
transformation. Also it is easy to see that $F$ is a composition
of a translation and a transformation from $GL_n(\HH)Sp(1)$. \qed

\section{The distribution $\left( \frac{\pt u}{{\partial\bar q_i \partial q_j}}\right)$.}
\setcounter{subsection}{1}
\setcounter{theorem}{0}
In this section we will define the  matrix valued measure
 $\left( \frac{\pt u}{{\partial\bar q_i \partial q_j}}\right)$
for any (finite) psh function $u$ on $\Omega$. This construction is a quaternionic version
of the well known analogous construction in the complex case (see \cite{lelong}, p.70).

\def\hn{{\cal H}_n}
\def\co{{\cal C}}
Let us denote by $\hn$ the (real) linear space  of $n\times n$ quaternionic
hyperhermitian  matrices. Let
$$\co:=\{\xi\in \hn |\, \xi\geq 0\}.$$
Thus $\co$ is a closed convex cone. On the space $\hn$ one has the bilinear
symmetric form $(\cdot, \cdot): \hn \times \hn \str \RR$ defined by
$$(A,B)=Re Tr(A\cdot B),$$
where for any $n\times n$ quaternionic matrix $X= (x_{ij})$, $Re
Tr (X):= \sum_{i=1}^n Re x_{ii}$. Note also that for any
quaternionic matrices $X$ and $ C$ with $C$ invertible, one has
$Re Tr(CXC^{-1})= Re Tr X$.

We easily have
\begin{claim}
(i) $(\cdot, \cdot)$ is a perfect pairing on $\hn$.

(ii) For any matrices $A,\, B\in \co$ one has $(A,\, B)\geq 0$.

(iii) The dual cone $\co ^{o}:=\{\xi\in \hn |(\xi,\eta)\geq 0 \,
 \forall \eta \in \co \} $ coincides with $\co$.
\end{claim}

\begin{definition}
One says that an $\hn$-valued distribution $\psi$ on the domain
$\Omega$ is {\itshape non-negative} ($\psi \geq 0$) if for any
smooth compactly supported function $f$ on $\Omega$  with values
in the cone $\co$ one has $\psi (f) \geq 0$. One can call such a
$\psi$ $\co$-valued.
\end{definition}
As in the usual scalar valued case one has the following result. (For
the scalar valued case see \cite{gelfand-shilov}.)
\begin{proposition}
Any $\co$-valued distribution on $\Omega$ is of zero order, i.e.
(non-negative) $\hn$-valued measure.
\end{proposition}
We also have the following result (which easily follows from the scalar valued
case).
\begin{lemma}
Any locally bounded sequence of $\hn$-valued measures on $\Omega$
has a weakly convergent subsequence.
\end{lemma}

Note also that for any $\hn$-valued distribution (resp. measure)
$\mu$ on $\Omega$ one can define in the obvious way its trace $Tr
\mu (=Re Tr \mu)$, which is a real valued distribution (resp.
measure).

\begin{proposition}
The sequence of $\co$-valued measures $\{\mu _j\}$ is (locally)
bounded iff the sequence of real valued measures $\{Tr \mu _j\}$
is (locally) bounded.
\end{proposition}
{\bf Proof.} It immediately follows form the fact that a subset
$X\subset \co (\subset \hn)$ is bounded iff the set $\{Tr x| \,
x\in X\}$ is bounded. \qed

Let us now define for any quaternionic psh function $u$ on
$\Omega$ the $\co$-valued measure $(\du)$ which has the usual
meaning for $C^2$- smooth function $u$. Let $u$ be an arbitrary
(finite) quaternionic psh function on $\Omega$. By \cite{alesker},
Subsection 2.1, $u$ is subharmonic. But every finite subharmonic
function is locally integrable (see e.g. \cite{ronkin}, Ch.1 \S
1.4). Hence we can define
$$u_{\eps}:= u\star \chi_{\eps},$$
where $\chi _{\eps}(z)= \frac{1}{\eps ^{4n}} \chi (\frac{z}{\eps})\geq 0$
is the usual smoothing kernel (like as in the complex situation, see e.g. \cite{hormander-CA} , p.45).
Then $u_{\eps}$ are $C^{\infty}$- smooth psh functions. Hence $(\due)\geq 0$ for all $\eps >0$.
\begin{proposition-definition}
For any quaternionic psh function $u$ on $\Omega$ the $\hn$- valued measures
$(\due)$ converge weakly to a non-negative $\hn$- valued measure as $\eps \str 0$.
This measure will be denoted by $(\du)$.
\end{proposition-definition}
It is easy to see that if $u\in C^2(\Omega)$ then the limit measure has
its usual meaning.

{\bf Proof.} First let us  show that the measures $\due$ are
locally bounded. By Proposition 3.1.5 it is sufficient to show
that their traces are locally bounded. But $Tr (\due)=\Delta
u_{\eps} \geq 0$. Let $K\subset \Omega$ be any compact subset. Let
$\gamma \geq 0$ be a smooth function with compact support on
$\Omega$ which is equal to 1 on $K$. Then
$$\int _K \Delta u_{\eps} d vol \leq \int _{\Omega} \Delta u_{\eps} \cdot \gamma d vol=
\int _{\Omega} u_{\eps}\cdot \Delta \gamma d vol \leq$$
$$||\Delta \gamma||_{C(\Omega)}\cdot \int _{supp \gamma} u_{\eps} d vol \leq
||\Delta \gamma ||_{C(\Omega)}\cdot \int _{supp (\gamma) +\eps D}|u| d vol =$$
$$ ||\Delta \gamma||_{C(\Omega)}\cdot ||u||_{L^1(supp (\gamma) +\eps D)},$$
where $D$ denotes the unit Euclidean ball in $\HH ^n$. This proves
the local boundedness of the sequence of measures $(\due)$. Hence
by Lemma 3.1.4 for any sequence $\{\eps _N\}\str 0$ the sequence
of measures $(\frac {\partial ^2 u_{\eps _N}} {\partial\bar q_i
 \partial q_j})$ has a weakly convergent subsequence. It
remains to show that the limit does not depend on the choice of
subsequence.

Fix an arbitrary $\phi \in C^{\infty}_{0}(\Omega)$. We easily get for any $i, \, j$
$$\int _{\Omega} \phi \cdot \due d vol=
 \int _{\Omega} u_{\eps} \cdot \frac{\partial ^2 \phi}{\partial\bar q_i \partial q_j} d vol
\str  \int _{\Omega} u \cdot \frac{\partial ^2 \phi}{\partial\bar q_i \partial q_j} d vol.$$
Namely for smooth $\phi$ the limit does not depend on the choice of subsequence.
This implies the statement. \qed

\begin{theorem}
(i) Let $\{u_N \}$ be a sequence of quaternionic psh functions on a domain
$\Omega \subset \HH^n$ which is uniformly bounded from above on every compact
subset of $\Omega$. Then either $u_N \str -\infty$ uniformly on compact subsets
of $\Omega$, or else there is a subsequence $\{u_{N_k}\}$ which converges in $L^1_{loc}(\Omega)$.
If $u_N \ne -\infty$ for all $N$ and $u_N$ converge in the sense of distributions to
a distribution $U$, then $U$ is defined by a psh function $u$ and $u_N \str u $
in $L^1_{loc} (\Omega)$.

(ii) Assume that  a sequence $\{u_N \}$ of quaternionic psh functions on $\Omega$ converges in
$L^1_{loc}(\Omega)$ to a quaternionic  psh function $u$. Then one has a weak convergence of measures
  $$(\dun)\str (\du).$$
\end{theorem}
To prove this theorem we will need a lemma which is a quaternionic analogue of the
corresponding complex result (see \cite{hormander-NC}, Theorem 4.1.7).
\begin{lemma}
Let $u$ be a function defined on a domain $\Omega \subset \HH ^n$.
Assume that $u_A(z) :=u(Az)$ is subharmonic in $\Omega _A :=\{z|
\, Az\in \Omega\}$ for every invertible linear quaternionic
transformation $A$ (i.e. $\forall A\in GL_n(\HH)$). Then $u$ is
quaternionic psh.
\end{lemma}
Assuming this lemma let us prove Theorem 3.1.7.

{\bf Proof} of Theorem 3.1.7. (i) This part of the theorem is known to be true if
one replaces in its statement the word "psh" by the word "subharmonic"
(see \cite{hormander-NC}, Theorem 3.2.12). In order to deduce part (i)
of the theorem from that result it remains to show that the limit function $u$
is psh (and not just subharmonic). But this immediately follows from Lemma 3.1.8.

(ii)  First note that the measures $(\dun),  \, N\geq 1$ are uniformly
locally bounded in $\Omega$. This is  proved exactly as in Proposition 3.1.6.

 Hence choosing a subsequence if necessary
we may assume that  this sequence of measures converges weakly to
an $\hn$-valued measure $(\nu _{\bar i j})$. We have to prove that
$\nu _{\bar i j}= \du$. To see it fix an arbitrary function
$\phi\in C^{\infty}_0(\Omega)$. Then
$$\int _{\Omega} \dun \cdot \phi d vol =
 \int _{\Omega} u_N \cdot \frac{\partial ^2 \phi}{\partial\bar q_i \partial q _j} d vol
\str \int _{\Omega} u \cdot \frac{\partial ^2 \phi}{\partial\bar q_i
\partial q _j} d vol= \int _{\Omega}\frac{\partial ^2 u}{\partial\bar q_i
\partial q _j}\phi \cdot dvol ,$$ where the first and the last equalities can be
easily deduced from the assumptions. The result follows. \qed

It remains to prove Lemma 3.1.8.

{\bf Proof} of Lemma 3.1.8. The proof is an easy modification of the proof of Theorem 4.1.7
in \cite{hormander-NC}. Fix $z\in \Omega$. The function
$u(z_1 +w_1, z_2 +\eps w_2, \dots , z_n + \eps w_n)$ is subharmonic in $w$ by
hypothesis, $0< \eps <1$. Hence
$$u(z) \leq \int _{|\zeta |=1} u(z_1 +r\zeta _1, z_2
 +r \eps \zeta _2, \dots ,z_n + r \eps \zeta _n) d \omega (\zeta),$$
where $\omega (\zeta)$ is the normalized Lebesgue measure on the unit sphere.
Since $u$ is upper semi-continuous and locally bounded above, the Fatou
lemma implies as $\eps \str 0$ that
$$u(z) \leq \int _{|\zeta|=1} u(z_1 +r \zeta _1, z_2, \dots , z_n) d \omega (\zeta).$$
The last inequality and Theorem 3.2.3 in \cite{hormander-NC} imply
that the function $z_1 \mapsto u(z_1,z_2, \dots ,z_n)$ is
subharmonic. The subharmonicity of the restrictions to other
quaternionic lines follows form the invariance under quaternionic
linear transformations. \qed


\section{The operator $\Phi(u)$.}
\setcounter{subsection}{1} \setcounter{theorem}{0} Following
Section 5 of the paper \cite{bedford-taylor} by Bedford and Taylor
we will define for any finite psh function $u$ an operator
$\Phi(u)$ which is essentially $(\det \pt u)^{\frac{1}{n}}$. It is
closely related to the operator $\det (\pt u)$, but it is defined
for arbitrary finite psh function $u$.

As in the previous section we will denote by $\co$ the  cone
of non-negative definite quaternionic hyperhermitian $n\times n$ matrices.
Consider the function
$$\Psi (\xi) = (\det (\xi))^{\frac{1}{n}}, \, \xi \in \co.$$
\begin{proposition}
The function $\Psi$ is a continuous, nonnegative, concave function which is
homogeneous of degree 1 on the cone $\co$.
\end{proposition}

{\bf Proof.} Concavity follows from Theorem 1.1.17 (ii) of \cite{alesker}.
The other properties are trivial. \qed

Let $\mu$ be a vector valued Borel measure on $\Omega \subset \HH^n$ with values in the cone $\co$.
 Let us define a nonnegative Borel
measure $\Psi(\mu)$ on $\Omega$ as follows. Choose a scalar valued
nonnegative Borel measure $\lambda$ on $\Omega$ so that $\mu$ is
absolutely continuous with respect to $\lambda$. Then by the
Radon-Nikodim theorem $d\mu= h \cdot d\lambda$ where $h$ is a
Borel measurable function on $\Omega$ with values in $\co$.
\begin{definition}
$\Psi(\mu):=\Psi(h) \lambda$.
\end{definition}
It is easy to see that this definition is independent of the choice of the measure $\lambda$. The following proposition is trivial.

\begin{proposition}
If $\mu$ and $\nu$ are Borel measures on $\Omega$ with values in $\co$ then

(1)$\Psi(\alpha \mu)= \alpha \Psi(\mu)$ if  $\alpha \geq 0$.

(2)If $\mu, \, \nu$ are mutually singular then $\Psi(\mu +\nu)=\Psi(\mu) +\Psi(\nu)$.

(3) $\Psi (\mu)$ is absolutely continuous with respect to $\mu$.

(4)$\Psi(t \mu +(1-t)\nu)\geq t\Psi(\mu) +(1-t)\Psi(\nu), \, 0<t<1$.
\end{proposition}

\begin{proposition}
If $\chi \geq 0$ is a continuous function with compact support then
$$\Psi(\mu \star \chi)\geq \Psi(\mu)\star \chi$$
on any compact set $\Omega '$ with $\Omega ' +support (\chi) \subset \Omega$.
\end{proposition}
{\bf Proof.} The proof is exactly the same as in the complex case
(see Proposition 5.4 in \cite{bedford-taylor}). It is essentially
based on  Proposition 4.1.3 and general measure theoretic
construction of Goffman and Serrin \cite{goffman-serrin}. We do
not reproduce it here. \qed

\begin{proposition}
Let $\mu_j$ be a sequence of Borel measures on $\Omega$ with
values in $\co$ which converges weakly to a Borel measure $\mu$.
Suppose also that Borel measures $\Psi(\mu _j)$ converge weakly.
Then
$$\Psi(\mu)\geq \underset{j\str \infty}{\lim} \Psi(\mu_j).$$
\end{proposition}
{\bf Proof.} Again the proof is exactly the same as in the complex case
(see Proposition 5.6 of \cite{bedford-taylor}). Note that in turn this is a special case of
\cite{goffman-serrin}, Theorem 3, page 165 with slight modifications. \qed

To define the operator $\Phi(u)$ for any psh function $u$ note that by Proposition-Definition 3.1.6 the
matrix of Borel measures $\left( \frac{\partial ^2 u}{\partial\bar q_i \partial q_j}\right)$
takes values in the cone $\co$.
\begin{definition}
$$\Phi(u) :=\Psi \left( \left( \frac{\partial ^2 u}{\partial\bar q_i \partial q_j}\right) \right).$$
\end{definition}
\begin{theorem}
Let $u,v,u_j$ be finite psh on $\Omega \subset \HH^n$. Then

(1) $\Phi(\alpha u)=\alpha \Phi(u),\, \alpha \geq 0,$

(2) $\Phi (tu+ (1-t)v)\geq t\Phi(u)+ (1-t)\Phi(v),\, 0<t<1,$

(3) If $\chi \geq 0$ is a continuous function with compact support then
$$\Phi(u\star \chi)\geq \Phi(u)\star \chi$$
on any open set $\Omega '$ with $\Omega ' + support(\chi)\subset \Omega$;

(4) if $u_j\str u$ as distributions on $\Omega$ and if  the
sequence of measures $\Phi(u_j)$ converges weakly then
$$\Phi(u)\geq \underset{j\str \infty}{\lim}\Phi(u_j).$$

(5) If $u_{\eps}=u\star \chi_{\eps}$ where $\chi
_{\eps}(z)=\frac{1}{\eps ^{4n}}\chi(\frac{z}{\eps})\geq 0$ is the
usual smoothing kernel (like in the complex situation, see e.g.
\cite{hormander-CA}, p.45) then
$$\underset{\eps \str 0}{\lim} \Phi(u_{\eps}) =\Phi(u).$$

(6) $\Phi(max\{u,v\})\geq min\{\Phi(u), \Phi(v)\}$.
\end{theorem}
{\bf Proof.} Assertions (1), (2) follow from Proposition 4.1.3.
Assertion (3) follows from Proposition 4.1.4. Let us prove (4). By
Theorem 3.1.7 $(\pt u_j)\str (\pt u)$ in the weak topology on the
space of $\hn$-valued measures on $\Omega$. This and Proposition
4.1.5 imply assertion (4). The assertions (5), (6) are proved
exactly as in the complex case, and we refer to the proof of
Theorem 5.7 in \cite{bedford-taylor}. \qed
\begin{theorem}
Let $n>1$. Let $u$ be a finite psh function on $\Omega$ such that the
regularizations of $u$, $u_{\eps} =u\star \chi_{\eps}$, have the
property that $\det (\pt u_{\eps})$ is a bounded family of Borel
measures on each compact subset of $\Omega$. Then

(1) $ \Phi(u)$ is absolutely continuous with respect to the Lebesgue measure,
 and if $\Phi(u) =g \cdot dvol$ then $g\in L^n_{loc}(\Omega)$, i.e. $g^n$ is locally
integrable;

(2) if $u$ is continuous and if $\det(\pt u)= f\cdot dvol +d\nu$
is the Lebesgue decomposition of the non-negative measure $\det(\pt u)$
into its absolutely continuous and singular parts then
$g^n\leq f$;

(3) if $ \frac{\partial ^2 u}{\partial\bar q_i \partial q_j}=
f_{\bar i j} +d\nu _{\bar i j}$ is the Lebesgue  decomposition of
the Borel measures $ \frac{\partial ^2 u}{\partial\bar q_i \partial
 q_j}$ then $g = (\det (f_{\bar i j}))^{\frac {1}{n}}$.
\end{theorem}
{\bf Proof.} The proof of this theorem is exactly as in the complex
case, and we refer to the proof of Theorem 5.8 in \cite{bedford-taylor}. \qed
\begin{remark}
The assumptions of Theorem 4.1.8 are satisfied when $u$ is a
{\itshape continuous} psh function.
\end{remark}

\section{On upper envelopes.}
\setcounter{subsection}{1} \setcounter{theorem}{0} Let
$\Omega\subset \HH^n$ be a domain. For functions $f\geq 0$ on
$\Omega$, and $\phi \in C(\bo)$ let us denote by
$$\bfm :=\{v \mbox{ is finite psh on }\Omega |\, \Phi (v)\geq f\cdot vol,
 \, \overline { \underset{q\str \zeta}{\lim}} v(q) \leq
\phi (\zeta)\, \forall \zeta \in \partial \Omega \}.$$
 The main result of this section is the following
result which is a quaternionic analogue of Theorem 6.2 from
\cite{bedford-taylor}.
\begin{theorem}
Let $\Omega$ be a strictly pseudoconvex bounded domain with smooth
boundary. Let $\phi \in C(\partial \Omega), \, f\in C(\bar
\Omega), \, f \geq 0$. Let
$$u(z):=\underset {v\in \bfm}{\sup  v(z)}.$$
Then $u\in C(\bar \Omega)$, $u$ is psh in $\Omega$, and
$u|_{\partial \Omega} \equiv \phi$. Moreover $u\in \bfm$.
\end{theorem}
The proof of this theorem closely follows \cite{bedford-taylor}
and \cite{bremermann}; we are going to present it. As in
\cite{bedford-taylor}, we will need two lemmas. Throughout this
subsection $\Omega$ will denote a strictly pseudoconvex bounded
domain in $\hn$.
\begin{lemma} Fix $\eps >0$. Then
for every $\zeta \in \bo$ there exists $\vz \in \bfm \cap C(\bar
\Omega) $ such that
$$\phi (\zeta) -\eps \leq \vz (\zeta) \leq \phi (\zeta).$$
\end{lemma}
\begin{lemma}
Fix $\eps >0$. Then for every $\zeta \in \bo$ there exists
$h_{\zeta} \in C(\bar \Omega)$ which is psh in $\Omega$ such that

1) $h_\zeta (z) \leq \phi (z)$ for all $z\in \bo$.

2) $h_\zeta (\zeta) \geq \phi (\zeta)-\eps$.
\end{lemma}
First let us show that Lemma 5.1.3 implies Lemma 5.1.2. Choose a
large constant $K\gg 0$ such that $\Phi (K|z|^2)= K\Phi
(|z|^2)>f$. Let $\tilde \phi (z):= \phi (z) -K(|z|^2-|\zeta|^2)$.
Let $h_\zeta$ be as in Lemma 5.1.3 applied for the function
$\tilde \phi$ instead of $\phi$. Then the function $\vz (z) :=
h_\zeta (z) +K(|z|^2-|\zeta|^2)$ satisfies Lemma 5.1.2.

 {\bf Proof} of Lemma 5.1.3. By assumption (at least locally in
a neighborhood of $\zeta$) $\Omega =\{F<0\}$ where
 $F$ is twice continuously differentiable function on $\HH ^n$
 which is strictly psh and $\nabla F|_{\bo} \ne 0$.
 Let $G(z):= F(z)- \delta |z- \zeta |^2$. For small $\delta$ the function $G$
 is psh in a neighborhood ${\cal O}$ of $\zeta$. Clearly
 $$G(\zeta) =0, \, G|_{(\bar \Omega -\{\zeta\})\cap {\cal O}}<0.$$
 Take small $\lam >0$ and consider
 $$F^*(z):= max\{G, -\lam\}.$$
 It is clear that

 1) $F^*$ is psh in $\Omega$ and continuous in $\bar \Omega$;

 2) $F^*(\zeta)=0$;

 3) $F^*|_{\bar \Omega -\zeta} <0.$
\newline
For given $\eps >0$ there exists a constant $C \gg 0$ such that
$$C\cdot F^*(z) +\phi(\zeta) <\phi(z) +\eps \mbox{ for all } z\in
\bo.$$ Let $ h_\zeta (z):=C\cdot F^*(z) +\phi(\zeta)-\eps$. Then
$$h_\zeta (z)<\phi(z)\mbox{ for all } z\in \bo,$$
$$h_\zeta (\zeta)=\phi(\zeta)-\eps.$$
This proves Lemma 5.1.3. \qed

{\bf Proof} of Theorem 5.1.1. Let us define the upper
regularization of $u$ in $\bar \Omega$ as usual:
$$u^*(z):= \underset{z'\str z}\limsup u(z').$$
It is easy to see that $u^*$ is psh in $\Omega$ (e.g. using Lemma
3.1.8 and the analogous classical result for subharmonic
functions, see \cite{ronkin} Ch. 1, \S 1.5) . Clearly $u\leq u^*$.

It follows from Lemma 5.1.2 that
$$u(z)\geq \phi(z) \forall z\in \bo.$$
In order to prove that $u$ coincides with $\phi$ on $\bo$ let us
prove the converse inequality (following \cite{bremermann}). Fix
$\zeta \in \bo$. As in the proof of Lemma 5.1.3 construct $F^*$.
Let $\fss :=-F^*$. Then

1) $\fss$ is super-harmonic;

2) $\fss (\zeta)=0$;

3) $\fss (z)>0$ for $z\in \bar \Omega -\zeta$.
\newline
In the classical potential theory $\fss$ is called {\itshape
barrier}, and hence the classical Dirichlet problem for harmonic
functions is solvable on $\Omega$. Hence there exists a harmonic
in $\Omega$ function $h\in C(\bar \Omega)$ such that $h|_{\bo}
\equiv \phi$. Since every function from $\bfm$ is subharmonic we
obtain that $u(z)\leq h(z)\, \forall z\in \bar \Omega.$ Since $h$
is continuous we get
$$u^*(z)\leq h(z)\mbox{ on } \bar \Omega,$$
$$\mbox{ hence }u^*(z)\leq \phi (z) \mbox{ on } \bo.$$
Finally we deduce
$$u(z)=u^*(z)=\phi (z), \, \forall z\in \bo.$$

By H. Cartan's theorem (see e.g. \cite{ronkin}) $u=u^*$ almost
everywhere in $\Omega$. Since $\bfm$ is closed under taking finite
maximums (Theorem 4.1.7 (6)) Choquet's lemma (see e.g.
\cite{ronkin}) implies that one can choose an increasing sequence
of functions $u_j\in \bfm$ which converges to $u$ almost
everywhere in $\Omega$. But then $u_j \str u^*$ in
$L^1_{loc}(\Omega)$. Hence by Theorem 3.1.7 (ii) $\pt u_j \str \pt
u^*$ weakly. Hence by Proposition 4.1.5
$$\Phi(u^*)\geq \underset{j\str\infty}{\lim} \Phi (u_j) \geq
\mu :=f\cdot vol =\mu.$$
Hence $u^*\in \bfm$. Since $u\leq u^*$ in $\bar \Omega$ we conclude
$$u\equiv u^*.$$
Hence $u$ is psh and $u\in \bfm$. Hence to finish the proof of
Theorem 5.1.1 it remains to prove the continuity of $u$ in $\bar
\Omega$.

First we will prove the following
\begin{claim}
$u$ is continuous at all points of the boundary $\bo$.
\end{claim}
{\bf Proof.} Fix any $\eps >0$ and any $\zeta \in \bo$. By Lemma
5.1.2 there exists a function $\vz \in \bfm \cap C(\bar \Omega)$
such that $\phi(\zeta)-\eps \leq \vz (\zeta)$. Since $\vz $ is
continuous, in a small neighborhood $U$ of $\zeta$ in $\bar \Omega$
$$\vz(z)>\phi(\zeta)-2\eps.$$
But $u(z)\geq \vz(z)$. Hence $u(z)> \phi(\zeta)-2\eps$ for $z\in
U$. Hence $\underset{z\str \zeta}{\liminf} u(z) \geq \phi(\zeta)$.
But $u$ is upper semi-continuous in $\clo$ (since $u\equiv u^*$),
hence $\underset{z\str \zeta}{\limsup} u(z)\leq
u(\zeta)=\phi(\zeta)$. Hence $u$ is continuous at $\zeta$. This
proves the claim. \qed

Now let us continue proving Theorem 5.1.1. Fix $\eps >0$. Let $\oe
>0$ be such that
$$\underset{ \underset{|z-z'|<\oe}{z,z'\in \bar \Omega}} {\sup}
|f(z)-f(z')|<\frac{\eps}{2} \mbox{ and } \underset{
\underset{|z-z'|<\oe}{z, z'\in \bar \Omega, dist(z,\bo)<3
\oe}}{\sup} |u(z)-u(z')|<\eps,$$ where $dist(z,\bo)$ denotes the
shortest distance from $z$ to $\bo$.
 Existence of such $\oe$ follows from Claim 5.1.4 and the continuity of $f$.
 Let $\tau
\in \HH^n$ be any vector with $|\tau|<\oe$ (where $|\cdot |$
denotes the norm of the vector). Let
$$v(z):=u(z+\tau)+\eps \cdot |z|^2- (L+1)\eps ,$$
where $L$ be any constant satisfying $L>|z|^2$ for all $z\in
\Omega$. Let
$$V(z,\tau):=\left\{ \begin{array}{c}
                           u(z), \mbox{ if } z\in \clo,\, z+\tau
                           \not\in \Omega \\
                           max\{u(z),v(z)\} , \mbox{ if }z\in \clo,\,z+\tau
                           \in \Omega
                     \end{array}
                 \right.
$$

\begin{lemma}
$V(z,\tau)\in \bfm$.
\end{lemma}
Let us postpone the proof of this lemma and let us finish the
proof of Theorem 5.1.1. Lemma 5.1.5 implies in particular that
$V(z,\tau)\leq u(z)$ for all $z\in \Omega $. Hence for any $z,\,
z+\tau \in \Omega$ such that $|\tau |<\oe$ we have
$$u(z+\tau)+\eps |z|^2-(L+1)\eps \leq u(z).$$
Hence for some constant $C$
$$ u(z+\tau)-u(z) <C\cdot \eps.$$
Replacing $\tau$ by $-\tau$ we get
$$|u(z+\tau)-u(z)| <C\cdot \eps.$$
Hence $u$ is continuous. \qed

 Thus it remains to prove Lemma 5.1.5.

{\bf Proof} of Lemma 5.1.5. Let us check all the conditions in the
definition of $\bfm$.
\begin{claim}
$V(z,\tau)\leq\phi(z)$ for all $z\in \bo$.
\end{claim}
{\bf Proof.} Indeed, if $z+\tau \not\in \Omega$ then
$V(z,\tau)=u(z)=\phi(z)$. If $z+\tau \in \Omega$ then either
$V(z,\tau)=u(z)=\phi(z)$ or
$$V(z,\tau)=u(z+\tau) +\eps \cdot|z|^2 -(L+1)\eps \leq
u(z+\tau)-\eps = $$ $$u(z)+ (u(z+\tau)-u(z)) -\eps \leq
u(z)=\phi(z).$$ \qed

Let us define the subset $\Gamma :=\{z\in \clo|\, z+\tau\in
\bo\}$. Note that for any point $x\in \Gamma$, $dist(x,\bo)\leq
\oe$. Let $A$ be the $\oe$-neighborhood of $\Gamma$. Then clearly
for all $x\in A$ one has $dist(x,\bo)\leq 2\oe$.
\begin{claim}
For all $z\in A$ one has $V(z,\tau)=u(z)$. Hence $V(z,\tau)$ is
upper semi-continuous in $\bo$.
\end{claim}
{\bf Proof.} Clearly it is sufficient to prove the first
statement. We have to check that for $z\in A$, $v(z)\leq u(z)$. We
have
$$v(z)=u(z+\tau) +\eps \cdot|z|^2 -(L+1)\eps \leq
u(z+\tau)-\eps\leq$$ $$u(z)+ (u(z+\tau)-u(z)) -\eps \leq u(z)$$
where the last inequality follows from the fact that
$dist(z,\bo)\leq \oe,\, |\tau|<\oe$ and the definition of $\oe$.
\qed

Since the maximum of two psh functions is psh we can easily get
from the last claim
\begin{claim}
$V(z,\tau)$ is psh in $\Omega$.
\end{claim}

To finish the  proof of  Lemma 5.1.5 it remains to prove
\begin{claim}
$\Phi(V(z,\tau)\geq \mu(=f\cdot dvol).$
\end{claim}
{\bf Proof.} Let us denote by $\ggo:= \Gamma \cap \Omega$. By
Claim 5.1.7 in a small neighborhood of $\ggo$, or if $z+\tau \not\in
\Omega$, we have
$$\Phi(V(z,\tau))=\Phi(u(z))\geq \mu.$$
Now it remains to consider domain $\{z\in \Omega|\, z+\tau \in
\Omega\}$. In this domain $V(z,\tau)=max\{u(z),v(z)\}$. Hence by
Theorem 4.1.7 (6) we get
$$\Phi(max\{u,v\})\geq min\{\Phi(u),\Phi(v)\}.$$
Since $\Phi(u)\geq \mu$ let us prove that $\Phi(v)\geq \mu$.
Indeed
$$\Phi(v(z))\geq \Phi(u(z+\tau)) +4\eps \geq f(z+\tau)+4\eps\geq
f(z).$$ \qed

Thus Lemma 5.1.5 and hence Theorem 5.1.1 are proved.
\subsection{Other Perron-Bremermann families.}
\def\Om{\Omega}
Let $\Omega$ be a domain in $\HH^n$. For brevity we will denote by
$P(\Omega)$ the class of psh functions in $\Omega$. Given $\phi\in
C(\bo)$ and a non-negative measure $\mu= f \cdot dvol$ on $\Om$,
we define three Perron-Bremermann families of subsolutions to the
Monge-Amp\`ere equation (the first one was defined earlier in
Subsection 5.1):

$\bfm:=\{v\in P(\Om)| \Phi(v)\geq\mu \mbox{ and } \underset{z\str
z_0}{\lim sup \,} \, v(z)\leq \phi (z_0),\mbox{ for all } z_0\in
\bo\}$;

$C\bfm:=\bfm \cap C(\bar \Om)$;

${\cal F} (\phi,\mu):=\{v\in P(\Om)\cap C(\bar \Omega)| det(\pt
v)\geq \mu \mbox{ and } v(z_0)\leq \phi (z_0) \mbox{ for all }
z_0\in \bo\}.$

If $\mu= f\cdot dvol,\, f\in L^n_{loc}(\Om)$ then let $\mu^n :=
f^n\cdot dvol$. If $v\in P(\Om)\cap C(\bar \Om)$ then by Theorem
4.1.8 (2) $$\Phi(v)^n\leq det(\pt v),$$ and consequently
$$C\bfm\subset {\cal F}(\phi,\mu^n).$$
By Theorem 4.1.8 (1) and Remark 4.1.9 if $u$ is continuous then
$\Phi(u)$ is absolutely continuous with respect to the Lebesgue
measure, $\Phi(u)= g\cdot dvol$, and $g\in L^n_{loc}(\Omega)$.
\begin{proposition}
Let $\Om$ be a bounded domain in $\HH^n$ and suppose that $u\in
P(\Om)\cap C(\bar \Om)$ satisfies $det(\pt u)=(\Phi(u))^n$. If
$C{\cal B} :=C{\cal B}(\phi,\Phi(u))$ and ${\cal F}:= {\cal
F}(\phi, det(\pt u)) $, where $\phi=u|_{\bo}$, then $\sup\{v|v\in
{\cal F}\}=\sup\{v|v\in C{\cal B}\}=u$.
\end{proposition}
{\bf Proof.} The remarks before this proposition and the
assumption imply that $C{\cal B}\subset{\cal F}$. Hence
$\sup\{v|v\in C{\cal B}\} \leq \sup\{v|v\in {\cal F}\}$. On the
other hand obviously $u\in C{\cal B}$. Hence $u\leq \sup\{v|v\in
C{\cal B}\}$. Thus it remains to show that $\sup\{v|v\in {\cal
F}\}\leq u$. Fix any $v\in {\cal F}$. By the minimum principle,
Theorem 2.1.3, $u-v$ attains its minimum on $\partial \Om$. But
since $u\geq v$ on $\bo$ we obtain that $u\geq v$ in $\bar \Om$.
The proposition is proved. \qed

\section{The Monge-Amp\`ere equation in the Euclidean ball.}
\setcounter{theorem}{0} In this section we prove the existence of
the solution of the Dirichlet problem for the  quaternionic
Monge-Amp\`ere equation for the unit Euclidean ball  assuming
sufficient regularity of the initial data $f,\, \phi$. First let
us introduce some notation.

Let $B$ denote the open unit ball in $\HH^n$,
$$B:=\{|q|< 1\}.$$
Let $\bar B$ be its closure. The main result of this section is as
follows:
\begin{theorem}
Let $f\in C^{\infty}(\bar B), \, f>0$. Let $\phi \in
C^{\infty}(\partial B)$. There exists unique psh function $u\in
C^{\infty}(\bar B)$ which is a solution of the Dirichlet problem
$$\det (\frac{\partial ^2 u}{\partial\bar q_i \partial q_j})=f,$$
$$u|_{\partial B}= \phi.$$
\end{theorem}

This (smooth) case will be used in the proof of the general case
(Theorem 0.1.3). The method of the proof of this case is a
modification of that of the paper by
Caffarelli-Kohn-Nirenberg-Spruck \cite{CKNS}.

 In this section we will denote
by $||g||_k$ the $C^k$- norm of a function $g$ in $\bar B$.

The proof uses the continuity method. In Subsection 6.1 we prove
the first order a priori estimates. In Subsection 6.2 we prove the
second order a priori estimates. In Subsection 6.3 we obtain
$C^{2,\alpha}$ a priori estimates as an easy consequence of the
results from Subsections 6.1 and 6.2 and a general result from
\cite{CKNS}. Then the higher smoothness results follow from these
by the standard regularity theory of elliptic equations of second
order (see e.g. \cite{gilbarg-trudinger}, \cite{krylov}).
\subsection{ First order estimates.}
\begin{proposition}
Assume that a psh function $u\in C^2(\bar B)$ satisfies the
quaternionic Monge-Amp\`ere equation with $f>0$ in $\bar B$. Then
$$||u||_1\leq C,$$
with a constant $C$ depending only on $||f||_1,\, ||\phi||_2$, and
$||f^{-1}||_0$.
\end{proposition}
{\bf Proof.} Let $L$ be the linearization of the operator
$v\mapsto log (det(\pt v))$ at $u$. Explicitly this operator can
be written
$$Lv= n f^{-1}\cdot det(\pt v,\pt u[n-1]).$$ Clearly $Lu=n$. Since $u$ is
strictly psh we have

{\bf Claim.}{\itshape The operator $L$ is elliptic.}
\newline
Let $D$ be a first order differential operator of the form
$D=\frac{\partial}{\partial x_i},$
where $x_i$ is one of the real coordinate axes in $\HH^n$. First
let us prove the following lemma.
\begin{lemma}
$$max_{\bar B}|Du|\leq max_{\partial B} |Du| +C,$$
where $C$ is a constant depending only on $||f||_1,\, ||\phi||_1$,
and $||f^{-1}||_0$ .
\end{lemma}
{\bf Proof.} We have $$L(Du)= nf^{-1}\cdot det(\pt(Du), \pt
u[n-1])=f^{-1}D(det(\pt u))=D(\log f).$$ Consider the function
$w(q):=\pm Du +\lam |q|^2$, with $\lam \gg 0$. Then we get
$$Lw= \pm D(\log f)+ \lam L(|q|^2).$$
But $$L(|q|^2)=8nf^{-1}\cdot det(I,\underset{n-1 \mbox{
times}}{\underbrace{\pt u,\dots, \pt u}})= 8n f^{-1}\cdot \sum
_{i=1}^n det(M_{ii}(\pt u)) ,$$ where $M_{ii}(A)$ denotes the
minor of a matrix $A$ obtained from $A$ by deleting the $i$-th row
and the $i$-th column.

{\bf Claim.} {\itshape Let $A$ be an invertible hyperhermitian
matrix of order $n$. For any $i,\, 1\leq i \leq n$,
$$(A^{-1})_{ii}=\frac{1}{det A} detM_{ii}(A).$$}

Thus using this claim we get
$$Lw = \pm D(\log f)+ 8n\lam  Tr(\pt u)^{-1}.$$
For any hyperhermitian positive definite $(n\times n)$ matrix $C$
one has $\frac{1}{n} Tr(C)\geq (det C)^{\frac{1}{n}}$. Hence we
get
$$Lw\geq \pm D(\log f)+ 8n^2 \lam \cdot (det (\pt u))^{-\frac{1}{n}}=\pm D(\log f)
+ 8n^2\lam f^{-\frac{1}{n}} .$$ Since $f\in C^1(\bar B)$ and $f$
is bounded from below by a positive constant, one can choose a
large $\lam$ such that the last expression will be positive. For
such a $\lam$ by the maximum principle the function $w$ achieves
its maximum on the boundary $\bb$. This proves Lemma 6.1.2.

Thus it remains to estimate the gradient $\nabla u$ on the
boundary $\bb$. First let $\tilde \phi$ denote any $C^2$- smooth
extension of $\phi$ inside the closed ball $\bar B$ such that its
$C^2$-norm can be estimated by the $C^2$-norm of $\phi$. Consider
the function $\tilde \phi +K(|q|^2-1)$ for large $K$. Let us
denote this extension again by $\phi$. Note that on the boundary
$\bb$ it coincides with our original $\phi$. Note also that for
large $K$ the function $\phi$ is psh and moreover
$$det(\pt \phi)\geq f= det(\pt u).$$
Hence by the minimum principle $\phi \leq u$ in $\bar B$. Next let $h$ be a
harmonic function in $B$ which extends $\phi$. Then $u\leq h$.
Hence on the boundary $|\nabla u|\leq max\{|\nabla h|,|\nabla
\phi|\}$. Thus Proposition 6.1.1 is proved. \qed

\subsection{Second order estimates.}
Let $D$ be any real first order differential operator with
constant coefficients which are not greater than one. First we
need the following result.
\begin{lemma}
For a constant $C$ depending only on $||f||_2, \, ||f^{-1}||_0$
$$max_{\bar B} D^2u \leq max _{\bb}D^2u +C.$$
\end{lemma}

{\bf Proof.} We have
$$D^2(\log f)=D^2(\log det(\pt u))=D(\frac{D(det\pt u)}{det \pt
u})=$$
$$ f^{-2}\left\{D^2(det \pt u)\cdot det(\pt u) -(D(det\pt u))^2\right\}=$$
$$f^{-2}\left\{ D[n\cdot det(\pt (Du),\pt u[n-1])]\cdot det(\pt u)-[n\cdot
det(\pt (Du), \pt u[n-1])]^2 \right\}=
$$
$$f^{-1} n\cdot det(\pt (D^2u), \pt u[n-1])+$$
$$ f^{-2} \left\{n(n-1)det (\pt (D u)[2],\pt u[n-2])\cdot det(\pt
u)-[n\cdot det(\pt (Du), \pt u[n-1])]^2 \right\}.$$ We need the
following

\begin{lemma} Let $A,\,B$ be hyperhermitian $(n\times n)$-matrices, $A>0$. Then
$$n(n-1)\cdot det(B[2], A[n-2])\cdot detA -(n\cdot det(B,
A[n-1]))^2 \leq 0.$$
\end{lemma}
Assuming this lemma let us finish the proof of Lemma 6.2.1. We get
$$D^2(\log f) \leq f^{-1} n\cdot det(\pt (D^2u), \pt u[n-1])=L(D^2 u).$$
Hence we have
$$L(D^2 u+\lam |q|^2)\geq D^2(\log f) +8n^2 \lam f^{-\frac{1}{n}},
$$
where we have used the lower estimate on $L(|q|^2)$ from the
previous subsection. For sufficiently large $\lam$ the last
expression is positive. Hence by the maximum principle
$$max_B (D^2 u+\lam |q|^2)\leq max_{\bb}(D^2 u+\lam |q|^2). $$
Thus Lemma 6.2.1 follows. \qed

{\bf Proof} of Lemma 6.2.2. The function $A\mapsto \log(detA)$ is
concave on the cone of positive definite hyperhermitian matrices
(see \cite{alesker}, Theorem 1.1.17 (i)).
 Hence $$\frac{d^2}{dt^2} (\log det(A+tB))|_{t=0} \leq 0.$$
Computing explicitly this derivative we obtain the lemma. \qed

Note now that in order to prove an estimate on the second
derivatives of $u$ it is sufficient to prove an upper estimate on
it. Indeed let $q_l=t+i\cdot x+j\cdot y+k\cdot z$ be one of the
quaternionic coordinates. Since $u$ is psh we have
$u_{tt}+u_{xx}+u_{yy}+u_{zz}\geq 0$. This and the upper estimates
on the second derivatives of the form $D^2u$ imply the lower
estimates on them. The estimates on the mixed derivatives also can
be obtained easily since
$$2u_{tx}=(\frac{\partial}{\partial t}+\frac{\partial}{\partial
x})^2u- u_{tt}-u_{xx}.$$

Hence we have to prove an upper estimate of $D^2u$ on $\bb$.
Let us introduce additional notation. Let $r(q)=|q|^2-1$. Then
$$B=\{r<0\}.$$ We will denote the quaternionic units as follows:
$$e_0=1,\, e_1=i,\, e_2=j,\, e_3=k.$$
 Fix a coordinate system $(q_1,\dots, q_n)$ on $\HH^n$; we
will write $q_i=\sum _{\eps =0}^3 e_\eps x_i^\eps$. Fix an
arbitrary point $P\in \bb$. We can choose such a coordinate system
near this point that the inner normal to $\bb$ at $P$ coincides
with the axis $x_n^0$. Also we will move the center of coordinates
to $P$, i.e. we will assume that $P$ coincides with $0$. Let us
denote the center of the ball $B$ by $R$.

First we have the following trivial estimates:
$$|u_{x_i^\eps,x_j^\delta}(P)|\leq C \mbox{ for } (i,\eps),\,
(j,\delta)\ne (n,0).$$ (Note that here we also use the first order
estimates of $u$ and $\phi$). Now let us prove the following
estimate:
\begin{lemma}
$$|u_{x_i^\eps,x_n^0}(P)|\leq C \mbox{ for } (i,\eps)\ne (n,0),$$
where $C$ depends only on $||f||_2,\, ||\phi||_3,\, ||f^{-1}||_0$.
\end{lemma}
{\bf Proof.} Clearly one can construct a vector field $T$ on
$\HH^n$ such that

1) $T(P)=\frac{\partial}{\partial x_i^\eps}$;

2) on the points of $\bb$, $T$ is parallel to $\bb$;

3) $T$ has the form
$$T=\frac{\partial}{\partial x_i^\eps}+a \cdot \frac{\partial}{\partial
x_n^0},$$ where the function $a$ is smooth  with estimates on the
derivatives depending only on $n$, and $a(P)=0$.

Consider the function
$$w(q):=\pm T(u-\phi)+ (u_{x_n^1} -\phi _{x _n^1})^2+(u_{x_n^2} -\phi _{x_n^2})^2+
(u_{x_n^3} -\phi _{x_n^3})^2- Ax_n^0 +B|q-R|^2.$$ We will show that
for $A,\,B$ sufficiently large
\begin{eqnarray*}
(a)Lw > 0;\\
(b) w|_{\bb} \leq 0.
\end{eqnarray*}
 If we will prove it then by the maximum principle
 $w\leq 0$ in $\bar B$. Hence
 $$|T(u-\phi)|\leq A x_n^0 \mbox{ in }\bar B.$$
 Hence at the point $P$,
 $|\frac{\partial}{\partial x_n^0} T(u-\phi)|\leq A.$
 This will finish the proof of Lemma 6.2.3.
 Thus let us check the conditions (a) and (b).
 By a
 straightforward computation
$$LT(u-\phi)= T(\log f) -(L(\phi_{x_i ^{\eps}})+a L(\phi _{x_n^0}))
+(u-\phi)_{x_n^0} La +$$
$$ n f^{-1} det\left((a_{\bar
i} \cdot (u-\phi)_{x_n^0,j}) + ((u-\phi)_{x_n^0,\bar i }\cdot
a_j), \pt u[n-1]\right).$$
 where we denote for brevity $g_i:=\frac{\partial g}{\partial
 q_i},\, g_{\bar i}:=\frac{\partial g}{\partial
 \bar q_i}$.
However $\frac{\partial g}{\partial x_n^0}=g_{\bar n}- e_1
\frac{\partial g}{\partial x_n^1}-e_2 \frac{\partial g}{\partial
x_n^2}-e_3 \frac{\partial g}{\partial x_n^3}$. Hence
$$LT(u-\phi)=T(\log f) -(L(\phi_{x_i ^{\eps}})+a L(\phi _{x_n^0}))
+(u-\phi)_{x_n^0} La +$$
$$nf^{-1}\left( det((a_{\bar
i}(u-\phi)_{\bar n,j})+ (a_{\bar i}(u-\phi)_{\bar n,j})^*,\pt
u[n-1])\right)-$$
$$nf^{-1}\sum_{l=1}^{3}det((a_{\bar i} e_l(u-\phi)_{x_n^l, j})+(a_{\bar
i} e_l(u-\phi)_{x_n^l, j})^*,\pt u[n-1]).$$
Using first and second order
estimates on $ a$ and $\phi$, and first order estimates on $f$ and
$u$ we get the following inequality:
$$|LT(u-\phi)|\leq C+ C nf^{-1} \det(I,\pt u[n-1]) +$$
$$nf^{-1}\left( det((a_{\bar i}u_{\bar n, j})+ (a_{\bar i}u_{\bar
n, j})^*,\pt u[n-1])\right)+$$
$$ nf^{-1}\sum_{l=1}^{3}det((a_{\bar
i} e_l (u-\phi)_{x_n^l, j})+(a_{\bar i} e_l (u-\phi)_{x_n^l,
j})^*,\pt u[n-1]).$$

We have the following linear algebraic identity:

{\bf Claim.}
$$det((a_{\bar i}u_{\bar
n,j})+ (a_{\bar i}u_{\bar n,j})^*,\pt u[n-1])= 2(Re \, a_{\bar n})
det(\pt u).$$

 It follows from Theorem 1.1.15 (i) of \cite{alesker} that for a fixed $n\times n$ positive definite
 hyperhermitian matrix $A$  the bilinear form
 $det(XX^*, A[n-1])$ is non-negative definite on the space of
 quaternionic $n$-columns. Hence we get
\begin{equation}
 |det(XY^*+YX^*, A[n-1])|\leq
 det(XX^*,A[n-1])+det(YY^*,A[n-1]).
\end{equation}

 Using this inequality and the last claim we obtain the following estimate:
 $$|LT(u-\phi)|\leq C+ C n f^{-1}\det(I,\pt u[n-1]) + 2n |Re \, a_{\bar n}|+$$
$$\sum _{l=1}^3 n f^{-1}\cdot \left( det(((u-\phi)_{x_n^l,\bar i}(u-\phi)_{x_n^l j}), \pt u[n-1]
) +det((a_{\bar i} a_{j}), \pt u[n-1])\right).$$
Using again the
first order estimates on  $a$,  we finally get
$$
| LT(u-\phi)|\leq C +C n f^{-1} det(I,\pt u[n-1]) +$$
$$ nf^{-1} \sum _{l=1}^3 det(((u-\phi)_{x_n^l,\bar i}(u-\phi)_{x_n^l j}), \pt
u[n-1]) ,\eqno{(2)}$$ but now the value of the constant $C$ might
be different from the previous one.

Now let us compute $L((u_{x_n^l}- \phi _{x_n^l} )^2)$. By  a
straightforward computation we have
$$
L((u_{x_n^l}- \phi _{x_n^l} )^2)= $$
$$2nf^{-1}
\det(((u-\phi)_{x_n^l,\bar i}(u-\phi)_{x_n^l, j
})+(u-\phi)_{x_n^l}\cdot ((u-\phi)_{x_n^l,\bar i j}), \pt
u[n-1])=$$
$$2nf^{-1}\det(((u-\phi)_{x_n^l,\bar i}(u-\phi)_{x_n^l,
 j}), \pt u[n-1])+$$
$$ 2(u-\phi)_{x_n^l}\cdot \left( (\log f)_{x_n^l}- n
f^{-1}\det(\phi _{x_n^l, \bar i j},\pt u [n-1]) \right).$$

Using this identity and (2)  we obtain:
$$Lw\geq -C+(8B-C)n f^{-1} \det(I,\pt u[n-1]) +$$
$$\sum _{l=1}^3 nf^{-1}\det(((u-\phi)_{x_n^l,\bar i}(u-\phi)_{x_n^l, j}), \pt
u[n-1])+$$
$$2\sum_{l=1}^3(u-\phi)_{x_n^l}\cdot \left( (\log
f)_{x_n^l}- n f^{-1}\det(\phi _{x_n^l, \bar i j},\pt u [n-1]) \right).$$
But the third summand is non-negative. Hence we get
$$Lw\geq -C+(8B-C)n f^{-1} \det(I,\pt u[n-1]) +$$
$$2\sum_{l=1}^3(u-\phi)_{x_n^l}\cdot \left( (\log
f)_{x_n^l}- n f^{-1}\det(\phi _{x_n^l, \bar i j},\pt u [n-1])
\right).$$ Using the first order estimates on $u$ and $f$  and
third order estimates on $\phi$ we finally obtain
$$Lw\geq -C'+ (8B-C')n f^{-1} \det(I,\pt u[n-1]).$$
As in the proof of Lemma 6.1.2
$$det(I,\pt u[n-1])\geq n f^{\frac{n-1}{n}}.$$
Thus for large $B$ we get
$$Lw\geq -C' +(8B-C')n^2 f^{-\frac{1}{n}}>0.$$

Thus the inequality (a) is proved. It remains to prove the
inequality (b), namely
$$(b)\, w|_{\bb} \leq 0$$
for large $A,\, B$. Note that $T(u-\phi)|_{\bb}\equiv 0$.
 Clearly it is sufficient to prove the inequality (b) only
near the point $P$. Since $u\equiv \phi $ on $\bb$, then using the
first order estimates on $u$ it is easy to see that for $l=1,2,3$
$$|u_{x_n^l}(q)-\phi _{x_n^l}(q)|< C|q|, \, q\in \bb.$$
But for $q\in \bb$ we have $|q|\leq K (x_n^0)^{\frac{1}{2}}.$
Hence $|u_{x_n^l}(q)-\phi _{x_n^l}(q)|^2< K' \cdot x_n ^0$. Thus
Lemma 6.2.3 is proved. \qed

Thus to obtain an estimate on all second order derivatives of $u$
it remains to prove
$$|u_{x_n^0,x_n^0}(P)|<C.$$
We have proven that $|u_{x_i^\eps, x_j^\delta}(P)|<C$ for
$(i,\eps),\, (j,\delta)\ne (n,0)$ and
$$|u_{x_n^l,x_n^0}(P)|<C \mbox{ for } l\ne 0. \eqno{(3)}$$
It suffices to show that $$|u_{\bar n n}(P)|<C.$$ However by $(3)$
it is  sufficient to show that for the $(n-1)\times (n-1)$- matrix
$$(u_{\bar \alpha, \beta}(P))_{\alpha, \beta <n} \geq c\cdot I
\eqno{(4)}$$ for some positive constant $c$. After subtracting a
linear functional we may assume that $\phi _{x_j^l}(P)=0$ for
$(j,l)\ne (n,0)$. In order to prove $(4)$ it is sufficient to
prove that
$$\sum_{\alpha, \beta <n} \bar \xi_{ \alpha} u_{\bar \alpha
\beta}(P) \xi_\beta \geq c|\xi|^2.$$ Let us prove it for
$\xi=(1,0,\dots,0)$. Namely $u_{\bar 11}\geq c$.

Let us write on the boundary $\bb$ the coordinate $x_n^0$ as a
function of other coordinates:
$$x_n^0=\rho((x_i^\eps)_{(i,\eps)\ne(n,0)}) .$$
\def\tu{\tilde u}
Let $\tu := u- \lam x_n^0$ with $\lam$ so chosen that
$$\Delta _1 \tu((x_i^\eps)_{(i,\eps)\ne(n,0)},\rho((x_i^\eps)_{(i,\eps)\ne(n,0)}))=0 \mbox{ at
} P ,$$
where $\Delta _1 =\sum _{\eps =0}^3 \frac{\partial ^2}{(\partial x_1^{\eps})^2}$.
 Since the first derivatives of $\rho$ vanish at $P$, the
last equality is equivalent to
$$\tu _{\bar 11}(P) +\tu_{x_n^0} \rho_{\bar 11}(P)=0.\eqno{(5)}$$
Consider the following Taylor decomposition:
$$\tu|_{\bb} =(\mbox{ quadratic terms in } x_i^\eps\ne
x_n^0)+(\mbox{ 3-order terms }) + O(|q|^4) = $$
$$E+F +O(\sum _{2\leq j \leq n}|q_j|^2) +O(|q|^4),$$
where
$$E:=(\mbox{ quadratic terms in } x_i^\eps\ne
x_n^0),$$
$$F:=(\mbox{ 3-order terms in }  x_1^{\eps}).$$
First let us consider the term $E$. We can estimate all the
monomials which do not contain $x_1^\eps$ by $C'\sum_{2\leq j\leq
n}|q_j|^2$. Thus
$$E= \sum_{\eps, \delta=0}^3 \sum_{\begin{array}{c}
                            j\ne 1\\
                             (j,\delta)\ne (n,0)
                             \end{array}}
a_{\eps, \delta,j}x_1^\eps x_j^\delta + Q(x_1^\eps) +O(\sum_{2\leq
j\leq n}|q_j|^2),$$ where $Q$ is a quadratic polynomial in
$x_1^\eps$ which satisfies $\Delta _1 Q =0$.

Now let us consider the expression $F$. It is well known (see
e.g. \cite{vilenkin}) that for any homogeneous polynomial $F$ of
degree 3 on a Euclidean space $\RR^N$ there exists a unique
decomposition $F(x)=F_0(x)+ l(x)\cdot |x|^2$, where $F_0$ is a
harmonic polynomial, and $l$ is a homogeneous polynomial of degree
1. Hence in our case ($N=4$) we can write
$$F= F_0+ (\sum_{\eps =0}^3 b^\eps x_1^\eps)|q_1|^2,$$
with $\Delta _1 F_0=0$.

On the boundary of the unit Euclidean ball $\bb$ we have
$$ 2x_n^0 = |q_1|^2 + \sum _{2\leq j \leq n-1} |q_j|^2 + \sum _{\delta
=1}^3 |x_n^\delta |^2 + O(|q|^3).$$ Thus
$$|q_1|^2=2x_n^0- (\sum _{2\leq j \leq n-1} |q_j|^2 + \sum _{\delta
=1}^3 |x_n^\delta |^2  )+ O(|q|^3).$$ Hence
$$F= F_0 + (\sum _{\eps =0}^3 b^\eps x_1^\eps) \left(2x_n^0 -(\sum _{2\leq j \leq n-1} |q_j|^2 + \sum
_{\delta =1}^3 |x_n^\delta |^2  )+O(|q|^3)\right)=$$
$$ F_0 +\sum _{\eps =0}^3 2b^\eps x_1^\eps
x_n^0 + O(\sum_{2\leq j \leq n}|q_j|^2 ) +O(|q|^4).$$

Thus we get an estimate
$$\tu|_{\bb}\leq \sum_{\eps ,\delta=0}^3 \sum_{\begin{array}{c}
                                 j\ne 1\\
                                   (j,\delta)\ne (n,0)
                                   \end{array}}
a_{\eps, \delta,j}x_1^\eps x_j^\delta + Q(x_1^\eps)+$$
$$(F_0 (x_1^\eps)+\sum _{\eps =0}^3 2b^\eps x_1^\eps x_n^0 ) +C \sum_{2\leq
j\leq n} |q_j|^2 +O(|q|^4),$$ where $\Delta _1 F_0=\Delta _1 Q=0$.
If we denote $G:=F_0+ Q$ then the last estimate can be rewritten
$$\tu|_{\bb}\leq \sum_{\eps,\delta=0}^3 \sum_{j\ne 1}
a_{\eps, \delta,j}x_1^\eps x_j^\delta +G+C \sum_{2\leq j\leq n}
|q_j|^2 +O(|q|^4).$$

\def\hu{\hat u}
 Let us define
$$\hu:= \tu -G.$$
Since $G$ depends only on $q_1$ and $\Delta _1 G=0$ then
$$\hu_{\bar i j}= \tu_{\bar i j}\mbox{ for } 1\leq i,j\leq n.$$
We have
$$\hu |_{\bb} \leq \sum_{\eps, \delta=0}^3 \sum_{j\ne 1}
a_{\eps, \delta,j}x_1^\eps x_j^\delta +C \sum_{2\leq j\leq n}
|q_j|^2 +O(|q|^4).$$

Now let us consider the following function
$$h:= -\alpha x_n^0 +\beta |q|^2 + \frac{1}{2D} \sum_{\eps,\delta}
\sum_{j\ne 1} |a_{\eps, \delta,j}x_1^\eps  +Dx_j^\delta|^2=$$
$$-\alpha x_n^0 +\beta |q|^2 +\sum_{\eps,\delta}
\sum_{j\ne 1}a_{\eps, \delta,j}x_1^\eps x_j^\delta +
D\sum_{\eps,\delta} \sum_{j\ne 1} |x_j^\delta|^2 +\theta ,$$ where
$\alpha,\, \beta$, and $D$ will be chosen later, and
$$\theta :=\frac{1}{2D}\sum_{\eps,\delta}
\sum_{j\ne 1}|a_{\eps, \delta,j}x_1^\eps|^2 \geq 0.$$ Hence
$$h|_{\bb} \geq -\alpha x_n^0 +\beta |q|^2 +\sum_{\eps,\delta}
\sum_{j\ne 1}a_{\eps, \delta,j}x_1^\eps x_j^\delta +4D
\sum_{j=2}^n |q_j|^2.$$ It is easy to see that for appropriate
choices of large $D$ and small $\alpha, \beta $ such that $-\alpha
x_n^0 +\beta |q|^2 \geq 0$ one can obtain that $h$ is psh and
$$h|_{\bb}\geq \hu|_{\bb}.$$
Now it is easy to see that the smallest eigenvalue of the matrix
$(h_{\bar i j})$ is equal to $4 \beta$. Clearly all the elements
of this matrix are bounded independently of small $\beta$; hence all
the other eigenvalues are bounded. Thus choosing sufficiently
small $\beta$ we may assume that
$$\det(h_{\bar i j}) <f \mbox{ in }  B.$$
Hence by the minimum principle
$$\hu \leq h \mbox{ in } \bar  B.$$ Since $h(P)=\hu(P)=0$ we
obtain
$$\hu _{x_n^0}(P) \leq h_{x_n^0}(P)=-\alpha.$$
It is easy to see that $\tu _{x_n^0}(P)=\hu _{x_n^0}(P)$.
Substituting this equality and the last inequality to $(5)$ we get
$$\tu _{\bar 11} (P)\geq \alpha \rho_{\bar 11}(P) =c>0.$$
But $u_{\bar 11} (P)=\tu _{\bar 11} (P)$. Thus the second order
estimate is proved. \qed

\subsection{$C^{2,\alpha}$- estimates.}
In this subsection we prove {\itshape a priori} $C^{2,\alpha}$-
estimates on solutions of the Dirichlet problem for the
quaternionic Monge-Amp\`ere equation. As previously we denote by
$u$ the solution of this problem. The main result of this
subsection is
\begin{theorem}
Let $\Omega$ be strictly pseudoconvex bounded domain in $\HH^n$
with smooth boundary. Let $u$ be a smooth psh solution of the
Dirichlet problem for the quaternionic Monge-Amp\`ere equation
with $f>0$, and $f,\, \phi$ be $C^{\infty}$-smooth. Then
$$|u|_{2+\alpha}\leq K \mbox{ for some } 0<\alpha <1,$$
where $K$ depends only on $\Omega$ and norms of $ f$ and $\phi$.
\end{theorem}

This theorem is an immediate consequence of the following general
result due to Caffarelli, Kohn, Nirenberg, and Spruck
\cite{CKNS} and the second order
estimates obtained in the previous subsection.
\begin{theorem}
Let $\Omega$ be a bounded domain in $\RR^N$ with the smooth
boundary $\bo$. Let $u$ be a smooth solution of the elliptic
equation
$$F(x,u,Du,D^2u)=0 \mbox{ in } \Omega,$$
$$u\equiv \phi \mbox{ on } \bo,$$
$\phi$ is smooth. Assume that $F$ is {\itshape concave } in the
second derivatives $u_{ij}$. Assume that $u$ satisfies an estimate
$$|u|_2\leq C'.$$
Then $$|u|_{2+\alpha}\leq K \mbox{ for some } 0<\alpha <1,$$ where
$K$ depends only on $\Omega, \,F,\,|\phi|_4, \,C'$.
\end{theorem}

Note that this theorem implies Theorem 6.3.1 if one takes
$F(x,u,Du,D^2 u)=\log(\det u_{\bar i j})-\log f$.

Thus Theorem 6.0.1 is proved as well.

\section{Proof of Theorem 0.1.3.}
\setcounter{subsection}{1} \setcounter{theorem}{0} In this section
we will finish the proof of our main result about existence of
solution of the Dirichlet problem for quaternionic Monge-Amp\`ere
equation (Theorem 0.1.3). Let the case $n=1$ is well understood in the literature, we will assume that $n>1$ throughout this section. 
But first we will need the following
result.
\begin{theorem}
Suppose $\Omega =B$ is the unit Euclidean ball in $\HH^n$,
$\phi \in C(\bb)$, $f\in C(\bar B), \, f\geq 0$. Let $d\mu =
f^{\frac{1}{n}} dvol$. Then the upper envelopes of the families
${\cal B}(\phi, \mu), \,C{\cal B}(\phi, \mu), {\cal F}(\phi, \mu ^n)$
coincide. If $u$ denotes the upper envelope, then
$u\in C(\bar B)$ and satisfies
$$\Phi(u)=f^{\frac{1}{n}} dvol \mbox{ in } B,$$
$$\det(\pt u)=f dvol  \mbox{ in } B,$$
$$u=\phi \mbox{ in } \bb.$$
\end{theorem}
{\bf Proof.}Remind that the Perron-Bremermann families from the
theorem were defined in Subsection 5.2. The argument follows very
closely the proof of Theorem 8.2 of \cite{bedford-taylor}. Choose
a sequence of functions $f_j
>0$ with $f_j\in C^\infty(\bar B)$ decreasing to $f$ uniformly on
$\bar B$. Choose also a sequence of functions $\phi _j\in
C^\infty(\bb)$ such that $\phi _j$ increases to $\phi$ uniformly
on $\bb$. By Theorem 6.0.1 there exist unique psh functions
 $u_j\in C(\bar B)$ which are solutions of the Dirichlet problem
$$\det(\pt u_j)=f_j  \mbox{ in } B, \, u_j=\phi_j \mbox{ in } \bb.$$
By the minimum principle (Theorem 2.1.3) the sequence $u_j$ is
increasing. We can choose positive numbers $\eta_j$ tending to
zero so that $\phi _j+\eta_j\geq \phi \mbox{ on } \bb$. Since
$$\det(\pt(u_k +\eps (|z|^2-1)))\geq \det(\pt u_k) +\eps^n
\det(\pt |z|^2)=f_k +\eps^n \det(\pt |z|^2),$$ and since $f_k\str
f$ uniformly we can choose positive numbers $\eps_j\str 0$ such
that
$$\det(\pt(u_k +\eps _j(|z|^2-1)))\geq \det(\pt u_j) \mbox{ for } k\geq j.$$
By the  minimum principle (Theorem 2.1.3) we get
$$u_k+\eps_j (|z|^2-1)\leq u_j(z)+\eta _j \mbox{ for } k\geq j, \, z\in \bar B.$$
But $u_j(z)\leq u_k(z)$. Hence $u_j\str u$ uniformly on $\bar B$.
By Theorem 2.1.2 $u$ is psh and $\det (\pt u_j)\str \det (\pt u)$
weakly. Hence $\det (\pt u)=f$. Further, by Theorem 4.1.7(4) $\Phi
(u)\geq f^{\frac{1}{n}} dvol $, and by Theorem 4.1.8 and Remark
4.1.9 $\Phi(u)^n\leq \det(\pt u)$. Hence
$\Phi(u)=f^{\frac{1}{n}}dvol$.

It follows from Proposition 5.2.1 that the upper envelopes of
$C{\cal B}(\phi, f \cdot dvol)$ and ${\cal F}(\phi, f^n\cdot
dvol)$ coincide with $u$. By Theorem 5.1.1 the upper envelopes of
${\cal B}(\phi, f\cdot dvol)$ and $C{\cal B}(\phi, f\cdot dvol)$
coincide. \qed

\begin{theorem}
Let $\Omega$ be a bounded open set in $\HH^n$. Let $\phi \in
C(\bo)$ and $d\mu =f^{\frac{1}{n}} \cdot dvol$ with $f\geq 0$ ,
$f\in C(\Omega)$. Suppose that

(i) ${\cal B}(\phi, \mu)$ is nonempty, and

(ii) the upper envelope $u=\sup \{v:v\in {\cal B}(\phi, \mu)\}$ is
continuous on $\bar \Omega$ with $u=\phi$ on $\bo$.
\newline
Then $u$ is psh and it is the solution to the Dirichlet problem
$$\det(\pt u)=f \mbox{ in } \Omega, \,u=\phi \mbox{ on } \bo.$$
Also $\Phi (u)=f^{\frac{1}{n}}\cdot dvol$.
\end{theorem}
{\bf Proof.} It should be checked only that $\det(\pt u)=f\cdot
dvol$ in $\Omega$. First let us show that $\det(\pt u)\geq f \cdot
dvol $ in $\Omega$. By Choquet's lemma there exists an increasing
sequence $u_j\in {\cal B}(\phi,\mu)$ which converges to $u$ almost
everywhere, and hence in $L^1_{loc}(\Omega)$. Then by Theorem
4.1.7(4) $\Phi(u)\geq f^{\frac{1}{n}}$. Let us write the Lebesgue
decomposition $$\det(\pt u)=\tilde f\cdot dvol +d\nu.$$ By Theorem
4.1.8(2)
$$f=(f^{\frac{1}{n}})^n\leq \tilde f.$$
Hence $\det (\pt u)\geq f\cdot dvol$.

 To prove the opposite inequality let us fix $z_0\in \Omega$, and choose
$\eps
>0$ so small that the closure of the ball $B(z_0,\eps)=\{|z-z_0|<\eps\}$ is contained in $\Omega$.
By the previous theorem there is a psh
function $v(z)\in C(\overline{B(z_0,\eps)})$ such that
$$v(z)=u(z) \mbox{ on } \partial B(z_0,\eps);$$
$$\Phi(v)=f^{\frac{1}{n}}\cdot dvol \mbox{ on } B(z_0,\eps);$$
$$\det(\pt v)=f \cdot dvol  \mbox{ on } B(z_0,\eps).$$
Since  $f=\det (\pt v)\leq \det (\pt u)$ on $B(z_0, \eps)$, by the
minimum principle we have $v\geq u$ in $\overline{B(z_0,\eps)}$.
Set $U(z)=v(z)$ if $z\in B(z_0, \eps)$, and $U(z)=u(z)$ if $z\in
\bar \Omega -B(z_0,\eps)$. Then clearly $U$ is continuous and psh,
and $U=\phi$ on $\bo$. We also have $\Phi (U)\geq f^{\frac{1}{n}}
\cdot dvol$. Therefore $U\in {\cal B} (\phi, f^{\frac{1}{n}}\cdot
dvol)$. Hence $U\leq u$. Hence $U\equiv u$. In particular in
$B(z_0,\eps)$ we have $\det(\pt u)=f$ and $\Phi
(u)=f^{\frac{1}{n}}\cdot dvol$. \qed

Finally let us prove Theorem 0.1.3.

{\bf Proof.} By Theorem 7.1.2 we have to  verify that ${\cal
B}(\phi, \mu)$ is not empty and its upper envelope $u\in C(\bar
\Omega)$, and $u=\phi$ on $\bo$. When $\Omega$ is strictly
pseudoconvex this is consequence of Theorem 5.1.1. Thus
$u=\sup\{v:v\in {\cal B}(\phi, \mu)\}$ is the solution of the
Dirichlet problem. \qed
\section{Quaternionic Levi form.}
\setcounter{subsection}{0} In this section we discuss some
additional properties of quaternionic strictly pseudoconvex
domains. In Subsection 8.1 we introduce a quaternionic version of
the Levi form of a domain with smooth boundary and prove that such
a domain is strictly pseudoconvex if and only if its Levi form is
positive definite. In Subsection 8.2 we consider some examples and
some other analogies with the real and complex cases. In
Subsection 8.3 we state some open questions.

\def\Om{\Omega}
\def\tz{T_{\bo, z}}
\def\htz{{}\!^{h}\tz}
\def\lz{L_{\bo, z}}
\subsection{The quaternionic Levi form.}
 \setcounter{theorem}{0} In this subsection we introduce the
quaternionic version of the Levi form. For the classical complex
case we refer to \cite{hormander-NC} and \cite{demailly}. The main
result of this subsection is Proposition 8.1.2.

 Let $\Om$ be a domain in $\HH^n$
with $C^2$-smooth boundary $\bo$. For any $z\in \bo$ let $T_{\bo,
z}$ denote the tangent space at $z$ to the boundary $\bo$. The
{\itshape quaternionic tangent space} to $\bo$ at $z$ is by
definition the maximal quaternionic subspace contained in $\tz$:
$$\htz :=\tz \cap \tz I \cap \tz J \cap \tz K.$$

Let $\rho\in C^2$ be a defining function of $\Om$, i.e.
$$\rho <0 \mbox{ on } \Om, \, \rho =0 \mbox{ and } d\rho \ne 0
\mbox{ on } \bo.$$

 The
{\itshape Levi form} $\lz$ on $\htz$ is defined as the restriction
of the hyperhermitian quadratic form $(\frac{\partial ^2 \rho (z)}
{\partial \bar q_i \partial q_j})$ to $\htz$ divided by $|\nabla
\rho (z)|$.
\begin{claim}
The Levi form does not depend on the choice of $\rho$.
\end{claim}

{\bf Proof.} Let $\rho '$ be another defining function of $\Om$.
Then in a small neighborhood of $z$ there exists a smooth function
$\alpha, \, \alpha (z)> 0$, such that $\rho'= \alpha \rho$. But
$$\frac{\partial ^2 (\alpha \rho) }{\partial\bar q_i \partial
q_j}=\alpha \frac{\partial ^2 \rho} {\partial\bar q_i \partial 
q_j}+ \frac{\partial \alpha}{\partial \bar q_i} \cdot
\frac{\partial \rho}{\partial q_j} + \frac{\partial \rho}{\partial
\bar q_i}\cdot \frac{\partial \alpha}{\partial q_j}+
\rho\frac{\partial ^2 \alpha } {\partial\bar q_i \partial 
q_j}.\eqno{(\star})$$ Now let us choose the coordinate system such
that $z$ is at the origin, and $\htz$ is spanned by the first
$n-1$ coordinates $q_1, \dots, q_{n-1}$. If we evaluate this
expression at $z$ we obtain for $i,j\leq n-1$:
$$\frac{\partial ^2 \rho' (z)}{\partial\bar q_i \partial q_j}=
\alpha \frac{\partial ^2 \rho (z)}{\partial\bar q_i \partial
q_j}.$$ \qed

\begin{proposition}
A $C^2$-smooth domain $\Om$ is strictly pseudoconvex iff the Levi
form if positive definite at each point $z\in \bo$.
\end{proposition}
{\bf Proof.} If $\Om$ is strictly pseudoconvex then there is
nothing to prove. Let us prove the opposite statement. Let us fix
a point $z\in \bo$ and let us assume that $\lz$ is positive
definite. Let us fix any defining function $\rho$ of $\Om$ in a
neighborhood of $z$. Let us also fix a coordinate system on
$\HH^n$ so that again $z$ is at the origin and $\htz$ is spanned
by the first $n-1$ coordinates $q_1, \dots, q_{n-1}$. From
$(\star)$ we obtain for any real valued smooth function $\alpha$:
$$\frac{\partial ^2 (\alpha \rho) }{\partial\bar q_i \partial
q_j}(z)=\alpha (z) \frac{\partial ^2 \rho(z)} {\partial\bar q_i
\partial q_j}+ \frac{\partial \alpha (z)}{\partial \bar q_i}
\cdot \frac{\partial \rho (z)}{\partial q_j} + \frac{\partial \rho
(z) }{\partial \bar q_i}\cdot \frac{\partial \alpha (z)}{\partial
q_j}.$$ Now let us choose $\alpha$ such that $\alpha
(q)=1+l(q_n)$, where $l$ is $\RR$-linear real valued functional
depending only on $q_n$. Then if either $i<n$ or $j<n$ we get
$$\frac{\partial ^2 (\alpha \rho)(z) }{\partial\bar q_i \partial
q_j}= \frac{\partial ^2 \rho(z)} {\partial\bar q_i
\partial q_j}.$$ For $i=j=n$ we get
$$\frac{\partial ^2 (\alpha \rho)(z) }{\partial\bar q_n \partial
q_n}=\frac{\partial ^2 \rho(z) }{\partial\bar q_n \partial q_n}+
\frac{\partial l(z)}{\partial \bar q_n} \cdot \frac{\partial \rho
(z)}{\partial q_n} + \frac{\partial \rho (z) }{\partial \bar
q_n}\cdot \frac{\partial l (z)}{\partial q_n}=$$
$$\Delta _n(\rho)+ 2Re (\frac{\partial l(z)}{\partial \bar q_n} \cdot \frac{\partial \rho
(z)}{\partial q_n}).$$ If we choose $l$ appropriately we can make
the last expression arbitrarily large, and then the matrix
$(\frac{\partial ^2 (\alpha \rho) }{\partial\bar q_i \partial
q_j})$ will be positive definite at $z$, and hence in some
neighborhood of $z$, and hence the function $\alpha \rho$ will be
strictly psh. But $\alpha \rho$ is also a defining functional of
$\Om$ near $z$. \qed

\subsection{Some examples.}
In this subsection we present a general construction of
quaternionic strictly pseudoconvex domains. It was suggested by M.
Gromov \cite{gromov} in analogy to the complex case. Then we
discuss some differences of the quaternionic situation with the
real and complex cases. This part depends very much on discussions
with M. Sodin.
\begin{definition}
Let $S$ be a real $3n$-dimensional linear subspace of $\HH^n$.
Then $S$ is called {\itshape totally real} if $$S\cap (S\cdot
I)\cap (S\cdot J)\cap (S\cdot K)=\{0\}.$$
\end{definition}
Note that a generic real $3n$-dimensional linear subspace is
totally real. Note also that $S$ is totally real if and only if
its orthogonal complement $S^\perp \subset (\HH^n)^*$ satisfies
$$S^\perp +I\cdot S^\perp +J\cdot S^\perp +K\cdot S^\perp
=(\HH^n)^*.$$
\begin{definition}
A smooth $3n$-dimensional submanifold  of $\HH^n$ is called
totally real if the tangent space at every point of it is totally
real.
\end{definition}
\begin{claim}
Let $M$ be a $3n$-dimensional totally real compact submanifold of
$\HH^n$. Let $\Omega :=M_{\eps}$ be the $\eps$-neighborhood of
$M$. Then for small $\eps >0$ the domain $\Omega$ is strictly
pseudoconvex.
\end{claim}

Now let us remind the following characterizations of convex (resp.
pseudoconvex) domains in $\RR^n$ (resp. $\CC^n$) (see e.g.
\cite{hormander-NC}).
\begin{claim}
Let $\Omega$ be a bounded domain in $\RR^n$ (resp. $\CC^n$). Then
$\Omega$ is convex (resp. pseudoconvex) if and only if the
function $x\mapsto -\log dist(x,\partial \Omega)$ is convex (resp.
plurisubharmonic).
\end{claim}

Unfortunately this criterion is {\itshape not} true in the
quaternionic situation already in $\HH^1$. Indeed by Proposition
8.1.2 any bounded domain with smooth boundary is strictly
pseudoconvex in the quaternionic sense. It is not difficult to
construct a domain $\Omega\subset \HH^1$ such that the function
$x\mapsto -\log dist(x,\partial \Omega)$ will be not subharmonic
(in the usual sense).

\subsection{Questions and comments.} We would like to state few
questions closely related to the material of this paper.

{\bf Question 1.} Find a geometric (or any other) interpretation
of the quaternionic Monge-Amp\`ere equation (or of an appropriate
modification of it).

Remind that the (modified) real Monge-Amp\`ere equations appear in
construction of convex hypersurfaces in $\RR^n$ with the
prescribed conditions on curvature. For this material we refer to
\cite{aubin}, \cite{pogorelov}. One of the main applications of
(modified) complex Monge-Amp\`ere equations is the construction of
K\"ahler metrics on complex manifolds. After the proof of the
Calabi-Yau theorem \cite{yau1}, \cite{yau2} and the Aubin-Yau
theorem \cite{aubin1}, \cite{yau2} they became the key tool in
complex differential geometry, see e.g. \cite{aubin},\cite{besse},
\cite{joyce-book} for further discussion.

{\bf Question 2.} (due to L. Polterovich.) Find a geometric
characterization of quaternionic strictly pseudoconvex domains.
(Note that we have not defined the notion of quaternionic
pseudoconvex domain in the non-strict sense.)


{\bf Question 3.} (due to G. Henkin.) This question is closely
related to the previous one. Let $\Omega \subset \HH^n$ be a
domain which admits an exhaustion by level sets of
plurisubharmonic function; in other words there exists a
plurisubharmonic function $h: \Omega \str \RR$ such that for any
number $c$ the set $\{h\leq c\}$ is compact. (Note that in the
classical complex situation this property is one of the equivalent
definitions of the pseudoconvex domain.) It was observed by G.
Henkin \cite{henkin} that if $h$ is strictly plurisubharmonic
Morse function then $\Omega$ admits a homotopy retraction onto a
compact subset of dimension at most $\frac{3}{4} \dim _{\RR}
\Omega =3n$ (indeed the Morse index of every critical point of
such a function is bounded from above by $3n$). This implies that
the boundary $\partial \Omega$ is connected provided $n>1$. These
properties are analogous to the corresponding properties of
pseudoconvex domains in the complex spaces (where the constant
$\frac{3}{4}$ is replaced by $\frac{1}{2}$). It would be of
interest to understand the relation between the class of domains
with this property and the class of strictly pseudoconvex domains
in the sense of this paper.

{\bf Question 4.} Generalize Theorem 6.0.1 on the existence of the
{\itshape regular} solution (under suitable  assumptions on
regularity of the initial data) to arbitrary strictly pseudoconvex
bounded domains with smooth boundary (and not only for the
Euclidean ball).

Note that the real analogue of this result was proved by
Caffarelli, Nirenberg, Spruck in \cite{CNS}, and the complex
analogue was proved by Caffarelli, Kohn, Nirenberg, Spruck in
\cite{CKNS} and Krylov in \cite{krylov2}.

\end{document}